\PassOptionsToPackage{colorlinks = true, citecolor = blue, urlcolor = blue}{hyperref}
\documentclass[11pt]{amsart}

\usepackage[T1]{fontenc}
\usepackage{amssymb}
\usepackage{mathtools}
\usepackage[numbers,sort&compress]{natbib}
\usepackage{doi}
\usepackage{tikz-cd}
\usetikzlibrary{decorations.pathreplacing,arrows.meta}
\usepackage{pgfplots}
\usepgfplotslibrary{groupplots}
\pgfplotsset{compat=1.18}
\pgfplotsset{
  manifold panel/.style={
    width=\linewidth,
    height=4.2cm,
    xmin=1.5, xmax=10.5,
    xtick={2,3,4,5,6,7,8,9,10},
    xtick pos=bottom,
    ytick pos=left,
    ytick style={draw=none},
    tick label style={font=\scriptsize},
    label style={font=\small},
    title style={font=\small},
    grid=major,
    ymajorgrids=true,
    xmajorgrids=false,
    grid style={draw=gray!30},
  },
}
\usepackage{todonotes}
\setuptodonotes{inline}
\usepackage{hyperref}
\usepackage{cleveref}
\usepackage{aliascnt}
\usepackage{enumitem}
\usepackage{pgfkeys}
\usepackage{booktabs}
\usepackage{float}

\usepackage[top=3.5cm, bottom=3.5cm, left=2.5cm, right=2.5cm]{geometry}

\theoremstyle{plain}
\newtheorem{theorem}{Theorem}[section]

\newaliascnt{lemma}{theorem}
\newtheorem{lemma}[lemma]{Lemma}
\aliascntresetthe{lemma}

\newaliascnt{proposition}{theorem}
\newtheorem{proposition}[proposition]{Proposition}
\aliascntresetthe{proposition}

\newaliascnt{corollary}{theorem}
\newtheorem{corollary}[corollary]{Corollary}
\aliascntresetthe{corollary}

\theoremstyle{definition}

\newaliascnt{definition}{theorem}
\newtheorem{definition}[definition]{Definition}
\aliascntresetthe{definition}

\newaliascnt{remark}{theorem}

\aliascntresetthe{remark}

\newaliascnt{example}{theorem}
\newtheorem{example}[example]{Example}
\aliascntresetthe{example}

\newaliascnt{experiment}{theorem}

\aliascntresetthe{experiment}

\crefname{theorem}{theorem}{theorems}
\Crefname{theorem}{Theorem}{Theorems}
\crefname{lemma}{lemma}{lemmas}
\Crefname{lemma}{Lemma}{Lemmas}
\crefname{proposition}{proposition}{propositions}
\Crefname{proposition}{Proposition}{Propositions}
\crefname{corollary}{corollary}{corollaries}
\Crefname{corollary}{Corollary}{Corollaries}
\crefname{definition}{definition}{definitions}
\Crefname{definition}{Definition}{Definitions}
\crefname{remark}{remark}{remarks}
\Crefname{remark}{Remark}{Remarks}
\crefname{example}{example}{examples}
\Crefname{example}{Example}{Examples}
\crefname{experiment}{experiment}{experiments}
\Crefname{experiment}{Experiment}{Experiments}

\crefname{appendix}{appendix}{appendices}
\Crefname{appendix}{Appendix}{Appendices}

\newcommand{\R}{\mathbb R}
\newcommand{\Z}{\mathbb Z}
\DeclareMathOperator{\conv}{conv}

\DeclareMathOperator{\supp}{supp}
\DeclareMathOperator{\diam}{diam}
\DeclareMathOperator{\vol}{Vol}
\DeclareMathOperator{\Ric}{Ric}
\DeclareMathOperator{\id}{id}
\DeclareMathOperator{\Cov}{Cov}

\newcommand{\Prok}{d_{\mathrm{Pr}}}
\newcommand{\Wass}{d_{\mathrm{W}}}
\newcommand{\dTr}{d_{\mathrm{L}}^{\,r}}

\newcommand{\card}{\#\,}

\begin{document}

\title{Wasserstein Stability of the Nerve Map}
\author{Kenneth McCabe and Jose A. Perea}
\thanks{This work was partially supported by the National Science Foundation through award \# DMS-2603389 and CAREER award \# DMS-2415445}

\begin{abstract}
\sloppypar
    Given a locally finite open cover of a paracompact space, a subordinate partition of unity defines a map into the geometric realization of the nerve. When the cover is good, the nerve theorem asserts that this map is a homotopy equivalence. In applied topology, one specializes this construction to covers of metric data by open balls centered at finitely many landmarks. The corresponding nerve map represents the data inside the \v{C}ech complex, enabling one to pull cohomological or spectral invariants back to the data. 
    Since point cloud data often comes equipped with weights encoding domain-specific importance, we introduce a weighted version of the nerve map and study its stability under perturbations of the landmarks and weights using metrics from optimal transport. 
    Our main result is a Lipschitz stability theorem for the weighted nerve map under a uniform nondegeneracy assumption. 
    We state the bound in terms of a scale-dependent optimal transport distance, and give an example where our estimate is attained. 
    Since this distance is dominated by a multiple of the $1$-Wasserstein distance, Wasserstein stability follows as a corollary. 
    We also give geometric conditions ensuring nondegeneracy, obtain explicit estimates of the Lipschitz constant for closed Riemannian manifolds, and sketch applications to Laplacian eigenvector cascading and Eilenberg--MacLane coordinates.
\end{abstract}

\maketitle

%\tableofcontents
\section{Introduction}
\label{Sec:Intro}
\sloppypar
Given a paracompact space $Z$ and a locally finite open cover $\mathcal U=\{U_\alpha\}_{\alpha\in\Lambda}$, a partition of unity $\{\phi_\alpha\}_{\alpha \in \Lambda}$ subordinate to $\mathcal U$ defines a   map from $Z$ into the geometric realization of the nerve of $\mathcal{U}$.
Explicitly, this map sends a point $z\in Z$ to the point with barycentric coordinates $(\phi_\alpha(z))_{\alpha\in\Lambda}$; we refer to it as the \emph{nerve map}. 
When $\mathcal U$ is good (i.e., every nonempty finite intersection of cover elements is contractible), the nerve theorem asserts that this map is a homotopy equivalence \cite[Corollary~4G.3]{hatcherAlgebraicTopology2002}.

In Topological Data Analysis (TDA), one often specializes this construction to  metric spaces and covers by finitely many open balls. 
Indeed, given a finite set of \emph{landmarks} $X$ in a metric space $M$ and a scale $r > 0$, the collection of open $r$-balls around points in $X$ cover its $r$-neighborhood $X^r\subseteq M$. 
A subordinate partition of unity then defines a map from $X^r$ into the geometric realization of the \v Cech complex of $X$ at scale $r$, enabling one to pull shape invariants, such as persistent cocycles or Laplacian eigenvectors, back to a neighborhood of the data.

A representative example of the utility of the nerve map is in constructing topological coordinates valued in \emph{Eilenberg--MacLane spaces}. For an Abelian group $G$ and $n\geq1$, an Eilenberg--MacLane space $K(G,n)$ represents cohomology in the sense that there is a natural bijection $[W,K(G,n)]\cong H^n(W;G)$ for any CW-complex $W$ \cite[Theorem~4.57]{hatcherAlgebraicTopology2002}. In particular, $K(\Z,1)\simeq S^1$, so composing the nerve map with a map representing an integer 1-dimensional  cohomology class produces circle-valued coordinates, as in the sparse circular coordinates algorithm of \cite{pereaSparseCircular2020}, following \cite{desilvaPersistentCohomology2011}.

For such coordinates to be meaningful in practice, one would like the nerve map to be stable under perturbations of the landmarks. Standard selection procedures, such as max-min sampling and random subsampling, can return different landmark sets $X$ and $Y$ after a change in initialization, or due to noise in the underlying dataset. As the nerve maps corresponding to different landmark sets take values in different \v{C}ech complexes, formulating stability of such maps requires a common metric space in which to compare them.\footnote{The \emph{Gromov--Hausdorff distance} \cite[Section~7.3]{burago2001course} overcomes this problem by infimizing the Hausdorff distance over all isometric embeddings into an arbitrary common metric space; one could make the case that a similar procedure could work here. Our inequalities would still hold after such infimation.} Following work of Adamaszek et al. \cite{AdamaszekAdamsFrick2018} on so-called \emph{metric thickenings}, we obtain a common space for both \v Cech complexes by identifying each point in the geometric realization with a finitely supported probability measure on the ambient metric space $M$, whose weights on the vertices are its barycentric coordinates. The nerve maps are then both valued in the space of probability measures on $M$, and thus can be compared using the \emph{$1$-Wasserstein distance} $\Wass$, a standard metric in the theory of optimal transport \cite{Villani2009OptimalTransportOldNew, santambrogio2015optimal}.

An additional observation also motivates the role of optimal transport in our work. Namely,  point cloud data often comes equipped with scalar quantities encoding domain-specific information, such as physical or chemical properties of atoms in biomolecular data \cite{mengWeightedPersistent2020a}, the luminosities of galaxies \cite{skibbaLuminosityweightedMarked2006}, or return intensities in LiDAR point clouds \cite{kashani2015lidar}. Such attributes can be used to assign nonnegative weights to points that encode their relative importance among the data. We therefore model weighted data as an empirical probability measure, following several influential TDA works \cite{chazalGeometricInference2011b,chazalRobustTopological2018,buchetEfficientRobust2016b}. To that end, we treat a weighted landmark set as a pair $(X,\omega_X)$, where $\omega_X$ is a probability measure whose support is $X$, and rescale the chosen partition of unity by $\omega_X$ to obtain a weighted version of the nerve map that recovers the classical map under uniform weights.

\subsection{Contributions}

Our main result, \Cref{Thm:Main-Theorem}, is a Lipschitz stability theorem for the weighted nerve map under a uniform nondegeneracy assumption. Let $M$ be a compact metric space of diameter $D$, let $X,Y\subseteq M$ be finite sets of landmarks with weights $\omega_X,\omega_Y$, and fix $r\in(0,D]$. Since the weighted nerve maps $\varphi_X^r$ and $\varphi_Y^r$ take values in a common Wasserstein space --- specifically, the measures on $M$ ---  we compare them uniformly on the intersection of their domains:
\[
d_\infty^{W_1}(\varphi_X^r,\varphi_Y^r)
:=\sup_{m\in X^r\cap Y^r}\Wass(\varphi_X^r(m),\varphi_Y^r(m)).
\]
On the weights, we use a scale-dependent refinement of $\Wass$. 
For each $m\in M$, we define a localized cost $c_m$ which ignores mass moved far away from $m$, since the nerve map at $m$ depends only on landmarks whose $r$-balls contain $m$. Taking the largest associated transportation cost $T_{c_m}$ gives
\[
\dTr(\omega_X,\omega_Y):=\sup_{m\in M}T_{c_m}(\omega_X,\omega_Y).
\]
We show that $\dTr$ is a metric  satisfying $\dTr\leq 3r\,\Wass$ and generating the same topology as $\Wass$ on the space of probability measures on $M$ (\Cref{Prop:Localized-Metric-Properties}).

\Cref{Thm:Main-Theorem} establishes Lipschitz bounds on $d_\infty^{W_1}(\varphi_X^r, \varphi_Y^r)$ in terms of $d_L^r (\omega_X, \omega_Y)$ and $\Wass(\omega_X, \omega_Y)$, which we refer to as the \emph{localized transport} and \emph{Wasserstein} estimates, respectively. The former estimate is attained in \Cref{Ex:Estimate-Sharp}, demonstrating the sharpness of our results. 
Each bound is obtained by relaxing a sharper quantity $\Delta^r_L(\omega_X,\omega_Y)$, which we call the \emph{pointwise estimate}, which is produced directly by the proof but fails to satisfy the triangle inequality, and thus is not a distance between the weights (see \Cref{Appendix:Fail}). Moreover, we give geometric conditions ensuring nondegeneracy (\Cref{Sec:Geo-Conditions}), and show that if  the reference measure is lower Ahlfors regular, then Voronoi-weighted nerve maps converge uniformly as the landmarks become denser (\Cref{Sec:Convergence}).

As a representative special case, we obtain the following result for Riemannian manifolds, proved in \Cref{Sec:Riemannian-Manifolds} as a special case of \Cref{Thm:Ricci}.

\begin{theorem}\label{Thm:Main-Intro}
Let $n\geq2$, let $M$ be a closed, connected $n$-dimensional Riemannian manifold of diameter $D$ with non-negative Ricci curvature, and let $\mu$ be the normalized volume measure. Let $X,Y\subseteq M$ be finite $\delta$-dense sets with $\mu$-measurable nearest-neighbor projections $p_X,p_Y$, and assign Voronoi weights $\omega_X=(p_X)_\#\mu$ and $\omega_Y=(p_Y)_\#\mu$. Then, for every $r\in(\delta,D]$,
\[
d_\infty^{W_1}(\varphi_X^r,\varphi_Y^r)
\leq \frac{(n+1)D^n}{(r-\delta)^{n+1}}\dTr(\omega_X,\omega_Y)
\leq \frac{3r(n+1)D^n}{(r-\delta)^{n+1}}\Wass(\omega_X,\omega_Y).
\]
\end{theorem}

Finally, we test our estimates on real and synthetic data (\Cref{Sec:Experiments}) and apply them to several topological coordinate pipelines (\Cref{Sec:Applications}). In our experiments, we observed that the pointwise estimate is quite close to the actual distance between nerve maps, typically within a factor of about $3$. Moreover, for large datasets, both the pointwise and localized transport estimates can be cheaper to compute exactly than the Wasserstein estimate (\Cref{Sec:Computational-Cost}).

\subsection{Related Work}

\subsubsection{Optimal Transport in TDA}
Several works in TDA have adopted a measure-theoretic view of point cloud data. Chazal et al. \cite{chazalGeometricInference2011b} replace the distance-to-set function by the \emph{distance-to-measure}, obtaining robustness to outliers and stability in the $2$-Wasserstein metric; subsequent variants have been studied extensively \cite{chazalRobustTopological2018,buchetEfficientRobust2016b,chazal2016rates,anai2019dtm,guibas2013witnessed,phillips2015geometric}. In a different direction, Blumberg et al. \cite{blumbergRobustStatisticsHypothesis2014} establish Gromov--Prokhorov stability for distributions of persistence barcodes from fixed-size samples. Chazal et al. \cite{chazalSubsamplingMethodsPersistent2015} obtain Wasserstein and Gromov--Wasserstein stability for expected persistence landscapes, while Hang et al. \cite{hangTopologicalStudy2018} prove Wasserstein stability for descriptors built from Fr\'echet functions. Following these works, we use optimal transport to compare weighted landmark sets and nerve maps.

\subsubsection{Metric Thickenings}
\label{Sec:Metric-Thickening-Work}
To make the stability of the nerve map precise, we require a metric on the geometric realization $|K|$ of a simplicial complex $K$. 
The usual weak topology need not be metrizable unless the complex is locally finite \cite[Remark~3.2]{AdamaszekAdamsFrick2018}, and familiar metrics need not account for the geometry of the vertex set. 
In contrast, the Wasserstein distance is geometry-aware in the sense that $\Wass(\delta_x,\delta_y)=d_M(x,y)$ for any $x,y$.  Identifying points of $|\check C(M;r)|$ with finitely supported probability measures and equipping them with $\Wass$ yields the \emph{\v Cech metric thickening} of \cite{AdamaszekAdamsFrick2018}.

Although the Wasserstein topology is coarser than the weak topology in general \cite{AdamaszekAdamsFrick2018}, the natural map $|\check C(M;r)|\to\check C^m(M;r)$ is always a weak homotopy equivalence \cite[Theorem~1]{gillespie_vietoris_2024}, and a homotopy equivalence when $M$ is a finite-dimensional Polish space \cite[Corollary~3.2]{AdamsKarassevVirk2025}; in particular, when $M$ is a Riemannian manifold. Whether the two spaces are homotopy equivalent for arbitrary metric spaces remains open.

\subsubsection{Wasserstein Transform}
Our work is also closely related to the Wasserstein Transform of Jin, M\'emoli, Smith, and Wan \cite{memoliWassersteinTransformArxiv2018}. A \emph{localization operator} assigns to each $m\in M$ a probability measure describing the data near $m$; the Wasserstein Transform then defines a pseudometric by comparing these measures in Wasserstein distance. 
Our weighted nerve map can be viewed as a localization operator with a compactly supported kernel, so our Wasserstein stability estimate is comparable to \cite[Lemma~48]{memoliWassersteinTransformArxiv2018}. However, by exploiting properties of our specific kernel, we remove the dependence on the diameter in that estimate, and we also obtain sharper bounds by refining the Wasserstein distance itself; see \Cref{Appendix:WT}.

\subsection{Outline}
\Cref{Sec:Preliminaries} gives the preliminary definitions we use throughout the paper. \Cref{Sec:Main} formulates and proves the main stability theorem, and gives several sharpness examples. \Cref{Sec:Geo-Conditions} gives geometric conditions ensuring nondegeneracy, explicit Lipschitz constants for Riemannian manifolds, and a convergence result. \Cref{Sec:Experiments} tests the sharpness and computational cost of our estimates on real and synthetic data. \Cref{Sec:Applications} elucidates the role of our estimates in the stability of several topological coordinate pipelines. \Cref{Sec:Discussion} summarizes the results and closes with directions for future work.

\section{Preliminaries}
\label{Sec:Preliminaries}

\subsection{Metric Measure Spaces}

Throughout the paper, we let $(M, d_M)$ be a compact metric space equipped with its Borel $\sigma$-algebra $\mathcal B(M)$, and write $D:=\diam(M)$. Throughout, nerve-map scales lie in $(0,D]$.

For $z\in M$ and $t>0$, let $B(z,t)$ (resp.\ $\overline B(z,t)$) denote the \emph{open} (resp.\ \emph{closed})  \emph{ball} of radius $t$ centered at $z$. For a subset $A\subseteq M$ and $\delta \geq 0$, the \emph{$\delta$-neighborhood of $A$} is
\[
A^{\delta} := \left\{ z \in M \; \Big| \; d_M (z,A) < \delta \right\}.
\]
The subset $A$ is said to be \emph{$\delta$-dense} if $d_M (z, A)\leq \delta$ for all $z\in M$.  The \emph{covering radius} of a finite $X\subseteq M$ is 
\[
\Cov (X) := \max\limits_{m\in M} d_M(m, X).
\]

A \emph{nearest-neighbor projection} of $M$ onto a finite subset $X$ is a surjection $p_X: M\to X$ with
\[
d_M (m,p_X(m))=\min_{x\in X} d_M (x,m) \quad \forall m \in M.
\]
By breaking ties according to a fixed total order on $X$, we can assume without loss of generality that nearest neighbor projections are Borel measurable. To see this, let $X= \{x_1, \dots, x_N\}$, and let $p_X$ be given by sending each $m \in M$ to its first nearest neighbor in this ordering. Then for every $i$ we have
\[
p_X^{-1}(\{x_i\}) = \bigcap_{j < i} \{ m \mid d_M(m, x_i) < d_M(m, x_j) \} \cap \bigcap_{ j > i} \{ m \mid d_M(m, x_i) \leq d_M(m, x_j) \},
\]
which is an intersection of Borel sets by continuity of the map $m \mapsto d_M(m, x_i) - d_M(m, x_j)$.

Let $\mathcal P(M)$ be the set of all Borel probability measures on $M$. For $\mu \in \mathcal P(M)$, let $f_\#\mu$ be the \emph{pushforward} of $\mu$ along a measurable function $f$, i.e., $(f_\#\mu)(A)=\mu (f^{-1}(A))$ for all Borel sets $A$.

\subsection{Metrics on Probability Measures}\label{Sec:Metrics-Prob-Measures}

We briefly review several standard metrics on probability measures.

\begin{definition}\label{Def:Coupling}
Let $\mu,\nu\in\mathcal P(M)$, and let $p_1,p_2:M\times M\to M$ denote the coordinate projections. A \emph{coupling} of $\mu$ and $\nu$ is a probability measure $\pi$ on $M\times M$ such that
\[
(p_1)_\#\pi=\mu,
\qquad
(p_2)_\#\pi=\nu.
\]
We denote the set of all couplings of $\mu$ and $\nu$ by $\Pi(\mu,\nu)$.
\end{definition}

Couplings are also referred to as \emph{transportation plans}, in view of the following construction.

\begin{definition}\label{Def:Transportation-Plan}
    Given probability measures $\mu, \nu \in \mathcal P(M)$ and a Borel measurable cost function $c : M \times M \to [0,\infty)$, the \emph{transportation problem} associated to $c$ is
\[
T_c (\mu, \nu) :=\inf_{\pi \in \Pi(\mu, \nu)} \int_{M \times M} c(x,y) \, d\pi(x,y).
\]
The \emph{$1$-Wasserstein distance} is the transportation problem $\Wass := T_{d_M}$.
\end{definition}

More generally, one has the $p$-Wasserstein distances $W_p$ for $p \in [1, \infty)$, so that $W_1 = \Wass$ \cite[Definition~6.1]{Villani2009OptimalTransportOldNew}. For $p = \infty$ one can let
\[
W_\infty(\mu,\nu):=\inf_{\pi\in\Pi(\mu,\nu)}\operatorname*{ess\,sup}_{(x,y)\sim\pi}d_M(x,y).
\]
We work with $p = 1$ throughout, since it is the only exponent that yields a Lipschitz stability result when used for both the nerve maps and the weights; see \Cref{Ex:p-Wasserstein}. For a comprehensive account of Wasserstein and related distances, we refer the reader to \cite{Villani2009OptimalTransportOldNew,santambrogio2015optimal}.

The Prokhorov distance provides a complementary way to compare probability measures.

\begin{definition}\label{Def:Prokhorov}
The \emph{Prokhorov distance} between $\mu, \nu \in \mathcal P(M)$ is
\[
d_{\mathrm{Pr}}(\mu,\nu)
:=
\inf \left\{
\epsilon > 0 :
\begin{array}{l}
\mu(A) \le \nu(A^\epsilon) + \epsilon,\\
\nu(A) \le \mu(A^\epsilon) + \epsilon
\end{array}
\text{ for every } A \in \mathcal B(M)
\right\}.
\]
\end{definition}

The Prokhorov and $1$-Wasserstein distances each metrize the weak convergence of measures on compact metric spaces \cite[Corollary~6.13 and following discussion]{Villani2009OptimalTransportOldNew}.

\subsection{Simplicial Complexes and Geometric Realization}
\label{Sec:SCpx-Geo-Real}

The \emph{geometric realization} of an abstract simplicial complex $K$ with vertex set $V$ is
\[
|K| := \left\{ \alpha : V \to [0,1]\; \Big| \; \supp \alpha \in K \text{ and } \sum_{v \in V} \alpha(v) = 1 \right\}.
\]

The \emph{\v Cech complex} of $X\subseteq M$ at scale $r \in (0, D]$ is the abstract simplicial complex
\[
\check C(X; r) := \left\{ \sigma \subseteq X \; \Big| \; 0 < \card \sigma < \infty \text{ and } \bigcap_{x \in \sigma} B_M(x,r) \neq \emptyset \right\}.
\]
Since the vertices of $\check C(X;r)$ lie in $M$, we may identify each $\alpha \in |\check C(X;r)|$ with the probability measure 
\[
\sum\limits_{x \in X} \alpha(x)\delta_x,
\]
whose support is a simplex of $\check{C}(X,r)$. This identifies $|\check C(X;r)|$ with the \emph{\v Cech metric thickening} of $X$ at scale $r$,
\[
\check C^m (X; r) := \{\mu \in\mathcal P(M): \supp\mu \in\check C(X;r)\}.
\]
Following \cite{AdamaszekAdamsFrick2018}, we regard $\check C^m (X;r)$ as a metric space by equipping it with the restriction of $\Wass$.

\subsection{Partitions of Unity}

Let $\mathcal U = \{ U_\alpha \}_{\alpha \in \Lambda}$ be a locally finite open cover of $M$. A \emph{partition of unity subordinate to $\mathcal U$} is a family of continuous functions $\{ \phi_\alpha : M \to [0,1] \}_{\alpha \in \Lambda}$ such that

\begin{itemize}
    \item[(i)] $\supp \phi_\alpha \subseteq \overline{U}_\alpha$ for every $\alpha \in \Lambda$;
    \item[(ii)] $\sum\limits_{\alpha \in \Lambda} \phi_\alpha(m) = 1$ for every $m \in M$.
\end{itemize}

For $t\in\R$, let $|t|_+:=\max\{0,t\}$, and for  $\rho>0$ define the triangular bump function
\[
k_\rho(m,x):=|\rho-d_M(m,x)|_+.
\]
Thus, $k_r(m, x) > 0$ if and only if $ d_M(m, x) < r$. Note that $k_r$ is $1$-Lipschitz in both arguments, and takes values in $[0,r]$.

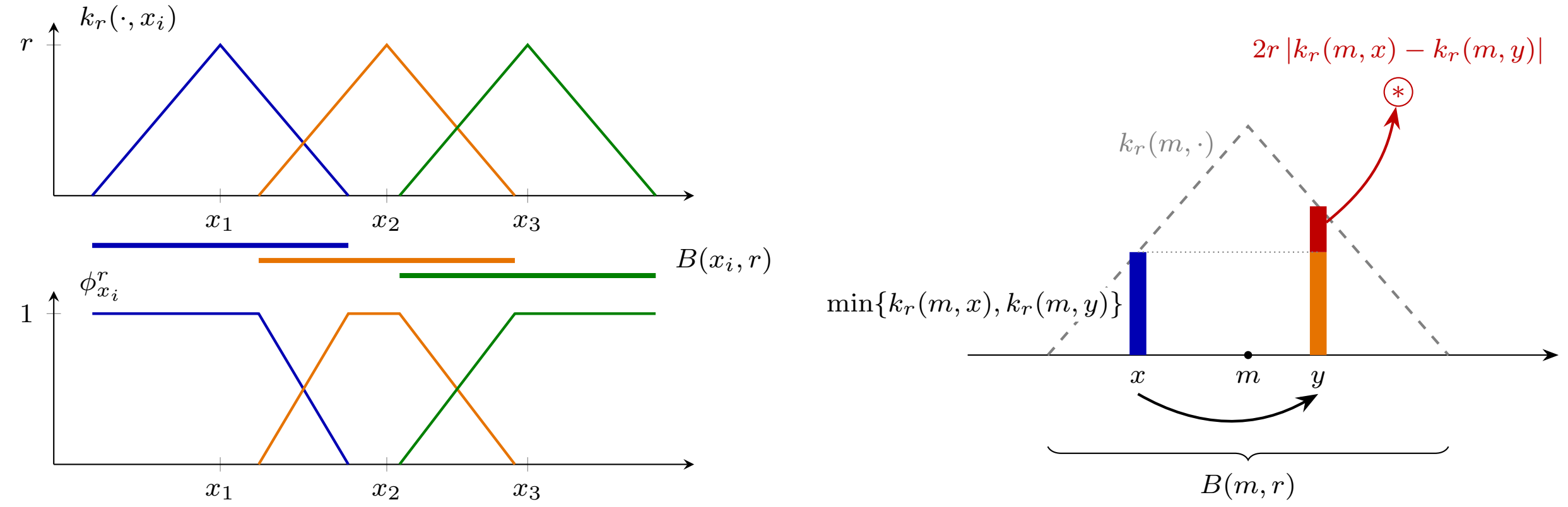
\begin{figure}[!tp]
\centering
% ---------------- Left: tents and partition of unity ----------------
\begin{minipage}[b]{.50\linewidth}
\centering
\begin{tikzpicture}[
  declare function={
    ka(\t) = max(0, 1 - abs(\t - 0));
    kb(\t) = max(0, 1 - abs(\t - 1.3));
    kc(\t) = max(0, 1 - abs(\t - 2.4));
    ks(\t) = ka(\t) + kb(\t) + kc(\t);
  }]
\colorlet{cA}{blue!70!black}
\colorlet{cB}{orange!90!black}
\colorlet{cC}{green!50!black}
\pgfplotsset{
  pou/.style={
    width=\linewidth, height=3.4cm,
    xmin=-1.3, xmax=3.7, ymin=0, ymax=1.15,
    axis x line=bottom, axis y line=left,
    xtick={0,1.3,2.4}, xticklabels={$x_1$,$x_2$,$x_3$},
    ytick={1}, yticklabels={$r$},
    tick label style={font=\scriptsize},
    clip=false, samples=241,
    every axis plot/.append style={thick},
  }
}
\begin{axis}[pou, name=top]
  \addplot[cA, domain=-1:1]   {ka(x)};
  \addplot[cB, domain=0.3:2.3] {kb(x)};
  \addplot[cC, domain=1.4:3.4] {kc(x)};
  % the r-balls, drawn as intervals below the axis
  \draw[cA, line width=1.5pt] (axis cs:-1,-0.33) -- (axis cs:1,-0.33);
  \draw[cB, line width=1.5pt] (axis cs:0.3,-0.43) -- (axis cs:2.3,-0.43);
  \draw[cC, line width=1.5pt] (axis cs:1.4,-0.53) -- (axis cs:3.4,-0.53);
  \node[font=\scriptsize, anchor=west] at (axis cs:3.45,-0.43) {$B(x_i,r)$};
  \node[font=\scriptsize, anchor=south west] at (axis cs:-1.2,1.0) {$k_r(\cdot,x_i)$};
\end{axis}
\begin{axis}[pou, at={(top.south)}, anchor=north, yshift=-1.0cm,
             ytick={1}, yticklabels={$1$}]
  \addplot[cA, domain=-0.999:1]  {ka(x)/ks(x)};
  \addplot[cB, domain=0.3:2.3]   {kb(x)/ks(x)};
  \addplot[cC, domain=1.4:3.399] {kc(x)/ks(x)};
  \node[font=\scriptsize, anchor=south west] at (axis cs:-1.2,1.0) {$\phi^r_{x_i}$};
\end{axis}
\end{tikzpicture}
\end{minipage}\hfill
% ---------------- Right: tent around m and the cost c_m ----------------
\begin{minipage}[b]{.48\linewidth}
\centering
\begin{tikzpicture}[x=2.1cm, y=2.4cm, >=Stealth]
\colorlet{cX}{blue!70!black}
\colorlet{cY}{orange!90!black}
\colorlet{cR}{red!75!black}
\def\xx{-0.55}\def\yy{0.35}
% axis and tent
\draw[->] (-1.4,0) -- (1.55,0);
\draw[gray, thick, dashed] (-1,0) -- (0,1) -- (1,0);
\node[gray, font=\scriptsize, anchor=east] at (-0.1,0.92) {$k_r(m,\cdot)$};
\fill (0,0) circle (1.2pt) node[below=1pt, font=\scriptsize] {$m$};
% landmarks x and y, with heights k_r(m,x)=0.45, k_r(m,y)=0.65
\draw[cX, line width=5pt] (\xx,0) -- (\xx,0.45);
\draw[cY, line width=5pt] (\yy,0) -- (\yy,0.45);
\draw[cR, line width=5pt] (\yy,0.45) -- (\yy,0.65);
\draw[gray, densely dotted] (\xx,0.45) -- (\yy,0.45);
\node[font=\scriptsize, below=1pt] at (\xx,0) {$x$};
\node[font=\scriptsize, anchor=east, fill=white, inner sep=1pt] at (\xx-0.05,0.22)
  {$\min\{k_r(m,x),k_r(m,y)\}$};
\node[font=\scriptsize, below=1pt] at (\yy,0) {$y$};
% common mass moved from x to y
\draw[->, thick] (\xx,-0.17) to[bend right=30] (\yy,-0.17);
% leftover mass sent to the auxiliary point
\node[circle, draw, cR, inner sep=1pt, font=\scriptsize] (star) at (0.75,1.15) {$\ast$};
\draw[->, thick, cR] (\yy+0.04,0.58) to[bend right=20] (star);
\node[cR, font=\scriptsize, anchor=south] at (star.north) {$2r\,|k_r(m,x)-k_r(m,y)|$};
% ball B(m,r)
\draw[decorate, decoration={brace, mirror, amplitude=4pt}]
  (-1,-0.4) -- node[below=4pt, font=\scriptsize] {$B(m,r)$} (1,-0.4);
\end{tikzpicture}
\end{minipage}
\caption{Left: $r$-balls, bump functions $k_r(\cdot,x_i)$, and the partition of unity $\phi^r_{x_i}$. Right: illustration of the cost function $c_m$. When a coupling pairs $x$ with $y$, the common mass (blue) is moved at cost $d_M(x,y)$ per unit, and the leftover mass (red) is sent to an auxiliary point $\ast$ at cost $2r$ per unit, as in \Cref{Lem:Cost-Pseudo}.}
\label{Fig:Tents-and-Cost}
\end{figure}

For $\mu\in\mathcal P(M)$ and $m\in M$, set
\[
S_\mu^r(m):=\int_M k_r(m,z)\,d\mu(z)
\]
and for $m\in(\supp\mu)^r$, define
\[
\varphi_\mu^r(m):=\frac{k_r(m,\cdot)}{S_\mu^r(m)}\,\mu\in\mathcal P(M).
\]
If $m \in (\supp \mu)^r$, then $k_r(m, \cdot)$ is positive on a neighborhood of some point of $\supp \mu$, so $S_\mu^r(m) > 0$. Conversely, if $m\notin(\operatorname{supp}\mu)^r$, then $k_r(m,z)=0$ for every $z\in\operatorname{supp}\mu$, and therefore $S_\mu^r(m)=0$. Thus $\varphi_\mu^r$ defines a map
\[
\varphi_\mu^r : (\supp \mu)^r \to \mathcal P(M),
\]
which we call a \emph{generalized nerve map}.

We are primarily interested in the case where $X \subseteq M$ is finite and $\mu = \omega_X$ satisfies $\supp \omega_X = X$. Here $(\supp \omega_X)^r = X^r$, and we write $S_X^r := S_{\omega_X}^r$ and $\varphi_X^r := \varphi_{\omega_X}^r$, so that
\[
\varphi_X^r(m) = \sum_{x\in X}\phi_{x}^r(m)\delta_x, \qquad \phi_x^r(m) := \frac{\omega_X(x)k_r(m,x)}{S_X^r(m)}, \qquad m \in X^r.
\]
The functions $\phi_x^r$ form a partition of unity subordinate to the cover of $X^r$ by the open $r$-balls around the landmarks $x\in X$, and when $\omega_X$ is uniform this is the standard such partition of unity.

\begin{center}
\begin{figure}[!tp]
    \includegraphics[width=.86\linewidth]{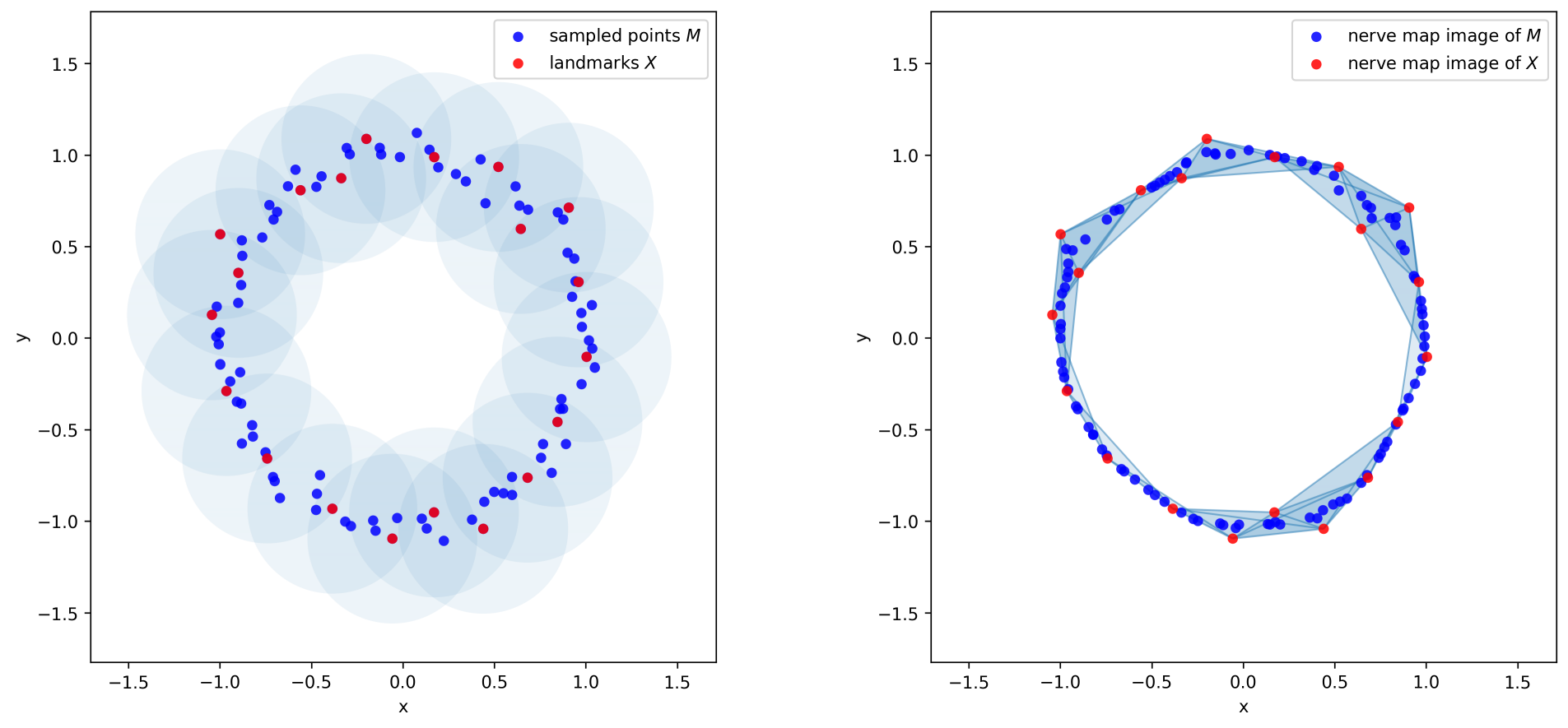}
    \caption{Left: $100$ points $M$ sampled from a noisy circle, with landmarks $X$ shown in red (equipped with uniform weights), and $r$-balls drawn in blue. Right: The image $\varphi_X^r(M)$ inside the \v Cech complex of $X$ at scale $r$.}
    \label{Fig:Weighted-Nerve-Map}
\end{figure}
\end{center}

Since $\phi_x^r(m)$ is positive only when $m$ lies in the $r$-ball around $x$, the support of $\varphi_X^r(m)$ spans a simplex of the \v Cech complex, and we obtain the \emph{weighted nerve map}
\[
\varphi_X^r : X^r \to \check C^m(X;r).
\]
See \Cref{Fig:Weighted-Nerve-Map,Fig:Tents-and-Cost} for examples.

\section{Stability of Weighted Nerve Maps}
\label{Sec:Main}

By viewing nerve maps as taking values in the common space $\mathcal P(M)$ of probability measures on $M$, we can compare them on the overlap of their domains.

Fix $r\in(0,D]$. For $\mu,\nu\in\mathcal P(M)$,
let
\[
\Omega_{\mu,\nu}^r
:=(\supp\mu)^r\cap(\supp\nu)^r,
\]
and define
\[
d_\infty^{W_1}(\varphi_\mu^r,\varphi_\nu^r)
:=\sup_{m\in\Omega_{\mu,\nu}^r}
\Wass(\varphi_\mu^r(m),\varphi_\nu^r(m)),
\]
with the convention $\sup\emptyset:=0$. For $\mu = \omega_X$ and $\nu = \omega_Y$, the supremum in the preceding definition is over $X^r \cap Y^r$.

We note that when the supports are allowed to vary arbitrarily, $d_\infty^{W_1}$ need not define a metric; indeed, the triangle inequality can fail, even when all domains have a nonempty intersection.

\begin{example}\label{Ex:Wass-Infty-Triangle}
    Let $M :=\{0, 1,3\} \subseteq \R$ with the Euclidean metric, and fix $r=2$. Define the measures 
    \[
    \mu := \frac 12 ( \delta_0 + \delta_3), \quad \nu := \frac 12( \delta_1+ \delta_3), \quad \lambda := \delta_3.
    \]
    Note that the nerve maps $\varphi_\mu, \varphi_\nu$ are defined on all of $M$, while $\varphi_\lambda$ is defined only at the point $3$. For $m \in \{0,1\}$, we have $\varphi_\mu^2(m)=\delta_0$ and $\varphi_\nu^2(m) = \delta_1$, while $\varphi_\mu^2(m)=\varphi_\nu^2(m) = \varphi_\lambda^2(m)=\delta_3$ at $m=3$. Therefore 
    \[
    d_\infty^{W_1}(\varphi^2_\mu, \varphi_\nu^2) = 1 > 0 = d_\infty^{W_1}(\varphi_\mu^2, \varphi_\lambda^2) + d_\infty^{W_1} (\varphi_\lambda^2, \varphi_\nu^2).
    \]
\end{example}

One can circumvent this pathology by fixing a common nonempty domain $\Omega$ over which to take the supremum. For example, if the nerve maps are defined on all of $M$, then one can take $\Omega = M$, so that $d_\infty^{W_1}$ is the usual uniform metric on the collection of nerve maps $M \to \mathcal P(M)$. For nerve maps corresponding to finite landmark sets with positive weights, a sufficient condition for this to hold is $\max\{\Cov(X),\Cov(Y)\}<r$. However, since none of our results require the triangle inequality for $d_\infty^{W_1}$, we retain the use of $\Omega_{\mu, \nu}^r$ in the supremum, which is the largest domain on which both nerve maps are defined and avoids imposing any additional coverage assumptions.

To motivate the metric we will use on weights, first observe that at a point $m \in M$, the nerve map depends only on the weights  of landmarks inside $B(m, r)$. We therefore define a localized cost function that assigns zero cost for a coupling to move mass between landmarks that lie outside $B(m,r)$. Indeed, for fixed $m\in M$ and $r\in(0,D]$, we define the cost function
\[
c_m (x,y) := \min \{ k_r(m,x),\, k_r(m,y) \}\, d_M (x,y) + 2r \,| k_r(m,x) - k_r(m,y)|.
\]
Both terms of $c_m$ have a geometric meaning: Suppose that a coupling of the weights pairs $x$ with $y$. The values $k_r(m, x)$ and $k_r(m, y)$ measure how much $x$ and $y$ contribute to the unnormalized nerve map at $m$. Moving the common mass $\min \{ k_r(m, x), k_r(m, y) \}$ from $x$ to $y$ costs the first term. The leftover $|k_r(m, x) - k_r(m,y )|$ units of mass must go somewhere in $B(m, r)$, a set of diameter at most $2r$, which costs at most the second term. This upper bound is exactly the cost of moving the leftover mass to an ``auxiliary'' point $\ast$ that is distance $2r$ away from every point of $B(m, r)$; see \Cref{Fig:Tents-and-Cost} and \Cref{Lem:Cost-Pseudo}. If $x$ and $y$ both lie outside $B(m, r)$, then $c_m(x,y)=0$, reflecting the fact that those landmarks do not contribute to the barycentric coordinates at $m$.

\begin{lemma}\label{Lem:Cost-Pseudo}
    For any fixed $r \in (0, D]$ and $m \in M$, the function $c_m$ is a pseudometric on $M$ with $c_m \leq 3r \, d_M.$
\end{lemma}

\begin{proof}
    Fix $m \in M$. Let $(Z_m, \hat d_m)$ be the metric space obtained from $(\overline B(m, r), d_M)$ by adjoining a point $\ast$ distance $2r$ from every point of $\overline B(m, r)$.

    Define the map $F_m :M  \to \mathcal P(Z_m)$ by
    \[
    F_m (x) := \frac 1r k_r(m, x)\delta_x + \left( 1- \frac 1r k_r(m, x) \right) \delta_\ast
    \]
    for $x\in B(m, r)$, and $F_m(x):= \delta_\ast$ otherwise.  Let $\Wass^{Z_m}$ denote the $1$-Wasserstein distance on $\mathcal P(Z_m)$, which is a metric since $Z_m$ is compact. It is straightforward to check that
    \[
    r \, \Wass^{Z_m} (F_m(x), F_m(y)) = c_m(x,y)
    \]
    for every $x, y\in M$.  It follows that $c_m$ is the pullback of a metric along $F_m$, and hence a pseudometric.

    Finally, since $k_r \leq r$ and $k_r(m, \cdot)$ is $1$-Lipschitz, we have
    \[
    c_m \leq r \, d_M + 2r \, d_M = 3r \, d_M. \qedhere
    \]
\end{proof}

We now use these cost functions to define a metric on $\mathcal P(M)$. For $\mu, \nu \in \mathcal P(M)$, let
\[
\dTr(\mu, \nu) := \sup_{m \in M} T_{c_m} (\mu, \nu).
\]

\begin{proposition}\label{Prop:Localized-Metric-Properties}
    For any fixed $r \in (0, D]$, the following hold:
    \begin{itemize}
        \item [(i)] $\dTr$ is a metric on $\mathcal P(M)$ with $\dTr \leq 3r \, \Wass$,
        \item [(ii)] $\dTr$ and $\Wass$ generate the same topology on $\mathcal P(M)$.
    \end{itemize}
\end{proposition}

\begin{proof}
    Fix $m \in M$. We first claim that $T_{c_m}$ is a pseudometric on $\mathcal P(M)$. Non-negativity and symmetry of $T_{c_m}$ are immediate from the corresponding properties of $c_m$, and the triangle inequality follows from a similar gluing argument as used in \cite[p.\ 94]{Villani2009OptimalTransportOldNew} to prove that $\Wass$ is a distance.  \Cref{Lem:Cost-Pseudo} gives $c_m\leq3r\,d_M\leq3rD$, and hence $T_{c_m}(\mu,\nu)\leq3rD$.  Thus, $\dTr$ is the supremum of a uniformly bounded family of pseudometrics, hence a pseudometric.

    It therefore remains to prove $\dTr$ separates points. Suppose $\dTr(\mu, \nu)=0$. For each $m \in M$, compactness of $M$ and continuity of $c_m$ give an optimal coupling $\pi_m \in \Pi(\mu, \nu)$ such that
    \[
    \int c_m \, d\pi_m = 0.
    \]
    Since $c_m \geq 0$, this implies that $c_m$ vanishes $\pi_m$-a.e.\ Moreover, $c_m(x,y)=0$ if and only if either $x=y$ or $x,y \notin B(m, r)$. Therefore, for every Borel $A\subseteq B(m, r)$, we have
    \[
    \mu(A) = \pi_m (A\times M) = \pi_m (A\times A) = \pi_m (M \times A) = \nu(A).
    \]
    Thus, $\mu$ and $\nu$ agree on every $r$-ball. Choosing a finite cover of $M$ by $r$-balls and refining it to a disjoint Borel partition yields $\mu = \nu$, completing the proof that $\dTr$ is a metric.

    For the second part of (i), fix $m \in M$ and a coupling $\pi$. Since $c_m \leq 3r \, d_M$, we have
    \[
    T_{c_m} \leq 3r \int d_M \, d\pi,
    \]
    so infimizing over couplings and taking the supremum over $m \in M$ gives $\dTr \leq 3r \, \Wass$.

    To prove (ii), the inequality $\dTr \leq 3r \Wass$ implies that the identity map
    \[
    \iota : (\mathcal P(M), \Wass) \rightarrow (\mathcal P(M), \dTr)
    \]
    is Lipschitz. Since $(\mathcal P(M), \Wass)$ is compact \cite[Remark~6.19]{Villani2009OptimalTransportOldNew} and $\dTr$ is a metric, it follows that $\iota$ is a continuous bijection from a compact space to a Hausdorff space, hence a homeomorphism.
\end{proof}

For $r \in (0, D]$ and $\mu, \nu \in \mathcal P(M)$, recall the notation $\Omega_{\mu, \nu}^r$ from above, and define
\[
\beta_\ast(\mu,\nu) := \inf_{m \in \Omega_{\mu, \nu}^r} \max \{ S_\mu^r(m), S_\nu^r(m) \}.
\]
For $m \in \Omega_{\mu, \nu}^r$, let
\[
D_m := \diam ( \supp \varphi_\mu^r(m) \cup \supp \varphi_\nu^r(m) ),
\]
and let
\[
\tilde c_m (x,y) := \min \{ k_r(m, x), k_r(m, y) \} d_M(x,y) + D_m | k_r(m, x) - k_r(m, y) |.
\]
We define the following function:
\[
\Delta_L^r (\mu, \nu) := \sup_{m \in \Omega_{\mu, \nu}^r} \frac{T_{\tilde c_m}(\mu, \nu)}{\max \{ S_\mu^r(m), S_\nu^r(m) \}}.
\]
We refer to this quantity as the \emph{pointwise estimate}. One can check that it is non-negative, symmetric, and reflexive, but since it fails the triangle inequality (see \Cref{Ex:Delta-Triangle}), it is not a metric in general. However, slight relaxations of $\Delta_L^r$ yield Lipschitz estimates in the metrics $\dTr$ and $\Wass$.

\begin{theorem}\label{Thm:Main-Theorem}
    Fix $r \in (0, D]$ and $\beta > 0$. Then for any $\mu, \nu \in \mathcal P(M)$ satisfying $\beta_\ast(\mu, \nu) \geq \beta$, we have
    \[
    d_\infty^{W_1}(\varphi_\mu^r, \varphi_\nu^r) \leq \Delta_L^r (\mu, \nu) \leq \frac{1}{\beta} \dTr (\mu, \nu) \leq \frac{3r}{\beta} \Wass(\mu, \nu).
    \]
\end{theorem}

\begin{proof}
    Write $\Omega := \Omega_{\mu, \nu}^r$. If $\Omega = \emptyset$, then the first two quantities vanish, so the result is immediate. Otherwise, assume $\Omega \neq \emptyset$.

    Fix $m \in \Omega$ and $\pi \in \Pi(\mu, \nu)$. By symmetry, assume without loss of generality that $S_\mu^r(m) \geq S_\nu^r(m)$, and put $S:= S_\mu^r(m)$, $\varphi_\mu := \varphi_\mu^r(m)$, $\varphi_\nu := \varphi_\nu^r(m)$, and $k(x) := k_r(m, x)$.

    Define a measure $\rho$ on $M \times M$ by
    \[
    d \rho (x,y) := \frac{ \min \{k(x), k(y)\}}{S} \, d\pi (x,y).
    \]
    Since $\pi$ has marginals $\mu,\nu$ and $\min\{k(x),k(y)\}\leq k(x),k(y)$, we have
    \[
    (p_1)_\#\rho\leq\frac{k}{S}\mu=\varphi_\mu,
    \qquad
    (p_2)_\#\rho\leq\frac{k}{S}\nu
    =\frac{S_\nu^r(m)}{S}\varphi_\nu\leq\varphi_\nu,
    \]
    where the last inequality uses $S\geq S_\nu^r(m)$. It follows that the residual measures
    \[
    \alpha_\mu := \varphi_\mu - (p_1)_\# \rho, \qquad  \alpha_\nu := \varphi_\nu - (p_2)_\# \rho
    \]
    are non-negative and have equal  mass $\eta:= 1 - \rho(M \times M)$. If $\eta >0$, then
    \[
    \tilde \pi := \rho + \frac{(\alpha_\mu \otimes \alpha_\nu)}{\eta}
    \]
    is a coupling of $\varphi_\mu$ and $\varphi_\nu$. If $\eta=0$, set $\tilde\pi=\rho$; the residual measures vanish, so $\rho$ already has the required marginals.

    Next, the residual measures are supported in $\supp \varphi_\mu \cup \supp \varphi_\nu$, a set of diameter at most $D_m$. Since $S= \int k(x) \, d\pi(x,y)$, we obtain
    \begin{align*}
        \Wass (\varphi_\mu, \varphi_\nu) &\leq \int d_M(x,y) \, d\rho(x,y) + D_m \eta \\
        &= \frac 1S \int \left( \min \{ k(x), k(y) \} d_M(x,y) + D_m ( k(x) - \min \{ k(x), k(y) \}) \right) \, d\pi(x,y) \\
        &\leq \frac 1S \int \tilde c_m(x,y) \, d\pi(x,y),
    \end{align*}
    where in the last line we used $k(x)-\min \{ k(x), k(y) \} \leq |k(x) - k(y) |$. Infimizing over $\pi$ and exchanging $\mu, \nu$ if necessary, we obtain
    \[
    \Wass (\varphi_\mu, \varphi_\nu) \leq \frac{T_{\tilde c_m} (\mu, \nu)}{\max \{S_\mu^r(m), S_\nu^r(m)\}}.
    \]
    Taking the supremum over $m \in \Omega$ proves the first inequality.

    For the remaining inequalities, since $\supp \varphi_\mu^r(m)$ and $\supp \varphi_\nu^r(m)$ lie in $\overline { B(m, r)}$ for every $m$, it follows that $D_m \leq 2r$, hence $\tilde c_m(x,y) \leq c_m (x,y)$. Moreover,
    \[
    \max \{ S_\mu^r(m), S_\nu^r(m) \} \geq \beta_\ast (\mu, \nu) \geq \beta.
    \]
    Consequently, \Cref{Prop:Localized-Metric-Properties}\,(i) gives
    \[
    \Delta_L^r(\mu, \nu)
    \leq \frac{1}{\beta} \sup_{m \in \Omega} T_{c_m}(\mu, \nu)
    \leq \frac{1}{\beta} \sup_{m \in M} T_{c_m}(\mu, \nu)
    = \frac{1}{\beta} \dTr(\mu, \nu) \leq \frac{3r}{\beta} \Wass(\mu, \nu). \qedhere
    \]
\end{proof}

Note that replacing $D_m$ by $2r$ and the pointwise normalization by $\beta$ gives the bound in $d_L^r$, and bounding $c_m \leq 3r\, d_M$ then gives the Wasserstein estimate.

We conclude this section with several examples. The first example shows that without a uniform lower bound on $\beta_\ast$, the ratio of $d_\infty^{W_1}(\varphi_X^r, \varphi_Y^r)$ to $\dTr(\omega_X, \omega_Y)$ can be arbitrarily large, even for two-point landmark sets in $\R$.

\begin{example}\label{Ex:Beta-Necessary}
Let $M = [-2,2] \subseteq \R$, $r = 1$, and for $\epsilon \in (0,1)$ let
\[
\omega_X := (1-\epsilon)\delta_{2} + \epsilon\,\delta_{1/2},
\qquad
\omega_Y := (1-\epsilon)\delta_{2} + \epsilon\,\delta_{-1/2}.
\]
At $m = 0 \in X^1 \cap Y^1$, only the landmarks $\pm 1/2$ lie within distance $1$, so $\varphi_X^1(0) = \delta_{1/2}$ and $\varphi_Y^1(0) = \delta_{-1/2}$, giving $d_\infty^{W_1}(\varphi_X^1, \varphi_Y^1) \geq 1$. On the other hand, $\Wass(\omega_X, \omega_Y) = \epsilon$, so $\dTr(\omega_X, \omega_Y) \leq 3\epsilon$ by \Cref{Prop:Localized-Metric-Properties}\,(i), and the ratio is at least $1/(3\epsilon) \to \infty$ as $\epsilon \to 0$. Here $\beta_\ast \leq S_X^1(0) = \epsilon/2 \to 0$.
\end{example}

The next example shows that the first two estimates of \Cref{Thm:Main-Theorem} are sharp, i.e., they can be equalities when $\beta = \beta_\ast$.

\begin{example}\label{Ex:Estimate-Sharp}
    Let $M$ be the vertex set of the star graph $K_{1,10}$ with the shortest path metric, let $m$ be the center vertex, and let $X:= \{ x_1, \dots, x_5\}$ and $Y:= \{ y_1, \dots, y_5\}$ be the leaves, equipped with uniform weights, and fix $r = 3/2$. Since $X^r \cap Y^r = \{ m\}$ and $k_r(m, z) = 1/2$ for every leaf $z$, we have $\beta_\ast = S_X^r(m) = S_Y^r(m) = 1/2$, and since $d_M(x_i, y_j) = 2$ for all $i,j$, $d_\infty^{W_1} (\varphi_X^r, \varphi_Y^r) = \Wass (\varphi_X^r(m), \varphi_Y^r(m)) = 2$.

    At $m$, we have $c_m (x_i, y_j) = \frac 12 \cdot 2 + 2r \cdot 0 = 1$ for all $i,j$, so $T_{c_m} (\omega_X, \omega_Y) = 1$. At a leaf, e.g., $x_1$, the only nonzero contribution to the cost comes from moving the $1/5$ mass at $x_1$, at cost $c_{x_1} (x_1, y_j) = 2r \cdot \frac 32 = \frac 92$, so $T_{c_{x_1}}(\omega_X, \omega_Y) = 9/10 < 1$, and by symmetry the same holds at every leaf. Hence $\dTr(\omega_X, \omega_Y) = 1$. At $m$, the residual term vanishes on $X\times Y$, so the two transportation costs coincide, and taking $\beta = \beta_\ast$, the first two inequalities of \Cref{Thm:Main-Theorem} collapse to
    \[
    2 = d_\infty^{W_1}(\varphi_X^r,\varphi_Y^r) \leq \Delta_L^r(\omega_X, \omega_Y) \leq \tfrac1\beta \dTr(\omega_X,\omega_Y) = 2.
    \]
\end{example}

Finally, we show that Lipschitz stability holds if and only if one uses the $1$-Wasserstein distance on the left-hand side, while any $q$-Wasserstein distance for $q \in [1, \infty]$ can be used on the right-hand side. For $p \in [1,\infty]$, define $d_\infty^p$ as $d_\infty^{W_1}$ with $W_p$ in place of $\Wass$.

\begin{example}\label{Ex:p-Wasserstein}
    Let $M = [-1/2, 3/2] \subseteq \R$ and $r = 1$, and for $0 < h < 1/2$ let $\omega_X = (\delta_0 + \delta_1)/2$ and $\omega_Y = (\delta_0 + \delta_{1-h})/2$. Then $X^r = Y^r = M$ and $\beta_* \geq \min_M S_X^r = \frac14$ for every $h$, while the coupling that fixes $0$ and moves $1$ to $1-h$ shows that $W_q(\omega_X, \omega_Y) = 2^{-1/q} h \leq h$ for all $q \in [1,\infty]$. Since $\varphi_X^1(0) = \delta_0$ and $\varphi_Y^1(0) = (\delta_0 + h \, \delta_{1-h})/(1+h)$,
    \[
    d_\infty^p(\varphi_X^1, \varphi_Y^1) \geq \Big( \frac{h}{1+h} \Big)^{1/p} (1-h) \geq \frac12 \Big( \frac{2h}{3} \Big)^{1/p}.
    \]
    For $p > 1$, the ratio of the two sides is thus at least $((2/3)^{1/p} h^{1/p-1})/2 \to \infty$ as $h \to 0$ with $r$ and $\beta_*$ fixed, so no estimate $d_\infty^p(\varphi_X^r, \varphi_Y^r) \leq C(r,\beta) \, W_q(\omega_X, \omega_Y)$ holds for any $q \in [1,\infty]$; for $p = \infty$ the left-hand side does not even tend to zero. For $p = 1$, by contrast, \Cref{Thm:Main-Theorem} and $W_1 \leq W_q$ \cite[Remark~6.6]{Villani2009OptimalTransportOldNew} give such an estimate for every $q$, and for $1<p<\infty$, since the nerve maps take values supported in sets of diameter at most $2r$, \cite[p.\ 179, Eq.~(5.1)]{santambrogio2015optimal} still yields a H\"older estimate with exponent $1/p$, which is sharp by this example.
\end{example}

\section{Geometric Conditions for Stability}
\label{Sec:Geo-Conditions}

In this section, we study the estimates of \Cref{Thm:Main-Theorem} quantitatively by comparing the weights with a reference measure $\mu$ satisfying a certain regularity property. We first obtain a uniform lower bound on $\beta_\ast$ when the weights are assumed to be close to $\mu$ (\Cref{Sec:Estimating-Beta}), and then make the resulting constants explicit for Riemannian manifolds (\Cref{Sec:Riemannian-Manifolds}). Finally, we show that the nerve maps converge as the weights converge to $\mu $ (\Cref{Sec:Convergence}).

\begin{definition}[Lower Ahlfors Regularity]\label{Def:Lower-Ahlfors-Regular}
We say that  $\mu \in \mathcal P(M)$ is \emph{lower Ahlfors $(c,q)$-regular} if there exist constants $c,q>0$ such that
\[
\mu(B(x,s))\ge cs^q,
\]
for every $x\in M$ and $s\in(0,D]$. When the constant $c$ is unimportant, we will also refer to $\mu$ as \emph{lower Ahlfors $q$-regular}.
\end{definition}

In this work, the main example of such a measure is the normalized volume measure on a closed, connected Riemannian $n$-manifold, where one can take $q = n$, and $c$ can be expressed in terms of the dimension, diameter, and Ricci curvature of $M$; see \Cref{Sec:Riemannian-Manifolds}. See \cite{heinonen2001lectures} for a general treatment of Ahlfors regular and related classes of measures on metric spaces.

\subsection{Comparing Weights with a Reference Measure}
\label{Sec:Estimating-Beta}

In order to make use of lower Ahlfors regularity to obtain explicit expressions for $\beta$, we first need to relate the normalization factor to ball masses. The following lemma is a direct application of the so-called ``layer cake'' representation for integrating non-negative measurable functions; see, e.g., \cite[Theorem~1.13]{liebAnalysis2001}.

\begin{lemma}\label{Lem:Layer-Cake}
    For any $m \in M$, $r \in (0, D]$, and $\omega \in \mathcal P(M)$, we have
    \[
    \int_M k_r(m,z)\,d\omega(z) = \int_0^r \omega(B(m, t))\, dt.
    \]
\end{lemma}

Thus, to find an admissible $\beta$ for \Cref{Thm:Main-Theorem}, it suffices to bound the ball masses of the weights from below and integrate.

\begin{proposition}\label{Prop:Ball-Mass-Stability}
    Let $X, Y \subseteq M$ be finite, let $\omega_X, \omega_Y \in \mathcal P(M)$ be probability measures whose supports are $X$ and $Y$, respectively, and fix $r \in (0, D]$. Suppose there exists a measurable function $g : [0,r] \to [0, \infty)$ such that, for every $m \in X^r \cap Y^r$,
    \[
        \min\{\omega_X(B(m,t)),\,\omega_Y(B(m,t))\}\geq g(t)
    \]
    for almost every $t \in [0,r]$. If $\beta := \int_0^r g(t) \, dt > 0$, then $\beta \leq \beta_\ast$, so the inequalities of \Cref{Thm:Main-Theorem} hold with this $\beta$.
\end{proposition}

\begin{proof}
    By \Cref{Lem:Layer-Cake}, $S_X^r(m) = \int_0^r \omega_X(B(m,t))\,dt \geq \beta$ for every $m \in X^r \cap Y^r$, and likewise for $S_Y^r(m)$. Hence $\beta \leq \beta_\ast$, and the conclusion follows from \Cref{Thm:Main-Theorem}.
\end{proof}

We now give two conditions under which \Cref{Prop:Ball-Mass-Stability} applies. In the first, we assume the weights are defined as pushforwards of the reference measure $\mu$. Geometrically, this amounts to assigning each landmark the $\mu$-mass of its Voronoi cell (\Cref{Fig:Voronoi-Weights}), up to measurable tie-breaking. In the second, we require only that the weights be Prokhorov-close to $\mu$. Here we use the Prokhorov distance as it gives direct control over ball masses, though one could use the inequality $\Prok^2 \leq \Wass$ \cite{gibbsChoosingBounding2002a} to obtain a similar (if messier) estimate for $\Wass$.

Let $X$, $Y$, $\omega_X$, and $\omega_Y$ be as in \Cref{Prop:Ball-Mass-Stability}, and fix   nearest neighbor projections $p_X : M \to X$ and $p_Y : M \to Y$, and a lower Ahlfors $(c,q)$-regular measure $\mu \in \mathcal P(M)$.

\begin{proposition}\label{Cor:Voronoi-Stability}
    Suppose $X$ and $Y$ are $\delta$-dense for some $\delta \geq 0$, and let $\omega_X := (p_X)_\#\mu$ and $\omega_Y := (p_Y)_\#\mu$. Fix $r \in (\delta, D]$, and set
    \[
    \beta := \frac{c}{q+1}(r-\delta)^{q+1} > 0.
    \]
    Then $\beta \leq \beta_\ast$, so the inequalities of \Cref{Thm:Main-Theorem} hold with this $\beta$.
\end{proposition}

\begin{proof}
    Fix $m \in X^r \cap Y^r$ and $t \in (\delta, r]$. Since $d_M(z, p_X(z)) \leq \delta$ for every $z \in M$, we have $B(m, t-\delta) \subseteq p_X^{-1}(B(m,t))$, so lower Ahlfors regularity gives
    \[
    \omega_X(B(m,t)) \geq \mu(B(m,t-\delta)) \geq c\,(t-\delta)^q,
    \]
    and same for $\omega_Y$. Using \Cref{Prop:Ball-Mass-Stability} with $g(t) := c\,|t-\delta|_+^q$ and integrating over $[0,r]$ gives $\beta$.
\end{proof}

\begin{proposition}\label{Cor:Prokhorov-Stability}
    Suppose that $\max\{\Prok(\omega_X,\mu), \Prok(\omega_Y,\mu)\} \leq \epsilon$ for some $\epsilon \geq 0$. Fix $r \in (\epsilon, D]$ with $c\,(r-\epsilon)^q > \epsilon$, and set
    \[
    \beta := \frac{c}{q+1}(r-\epsilon)^{q+1} - \epsilon\,(r-\epsilon) + \frac{q}{q+1}\,c^{-1/q}\,\epsilon^{(q+1)/q} > 0.
    \]
    Then $\beta \leq \beta_\ast$, so the inequalities of \Cref{Thm:Main-Theorem} hold with this $\beta$.
\end{proposition}

\begin{proof}
    Fix $m \in X^r \cap Y^r$ and $t \in (\epsilon, r]$. For every $\epsilon' > \epsilon$, by definition we have $\mu(B(m,t-\epsilon)) \leq \omega_X(B(m,t-\epsilon+\epsilon')) + \epsilon'$. Let $\epsilon' \downarrow \epsilon$, and note that $\omega_X(\bar B(m,t)) = \omega_X(B(m,t))$ for all but finitely many $t$, since $\omega_X$ is finitely supported. This gives $\omega_X(B(m,t)) \geq c\,(t-\epsilon)^q - \epsilon$ for almost every $t$, and likewise for $\omega_Y$. Define
    \[
    g(t):=\max\{c\,|t-\epsilon|_+^q-\epsilon,0\},\qquad 0\leq t\leq r.
    \]
    A direct calculation gives $\int_0^r g(t)\,dt = \beta$, and the conclusion follows from \Cref{Prop:Ball-Mass-Stability}.
\end{proof}

\subsection{Estimates for Riemannian Manifolds}
\label{Sec:Riemannian-Manifolds}

We now derive explicit Lipschitz constants for the estimates of \Cref{Thm:Main-Theorem} when $M$ is a Riemannian manifold. More specifically, let $M$ be a closed, connected Riemannian manifold of dimension $n \geq 2$, equipped with the geodesic distance, and let $\mu$ denote the normalized volume measure.

\begin{proposition}[\cite{MSEPost}]
\label{Prop:Ahlfors-Riemannian}
Let $\Ric_M$ denote the Ricci curvature of $M$, and choose $K \in \mathbb R$ such that $\Ric_M\ge (n-1)K$.\footnote{Such a $K$ exists since the Ricci tensor is continuous and the unit tangent bundle is compact \cite[Proposition 2.9 and Lemma 7.15]{lee_introduction_2018}.} Let $v_n$ be the volume of the unit ball in $\mathbb R^n$, and let $\vol_K(r)$ be the volume of a geodesic ball of radius $r$ in the complete, simply-connected $n$-dimensional Riemannian manifold of constant sectional curvature $K$. Then $\mu$ is lower Ahlfors $(c,n)$-regular, with $c = D^{-n}$ if $K \geq 0$ and $c = v_n/\vol_K(D)$ if $K < 0$.
\end{proposition}

\begin{proof}[Proof sketch]
    Following Yang \cite{MSEPost}, the Bishop--Gromov inequality \cite[Lemma~7.1.4]{petersenRiemannianGeometry2016}, applied with $R \downarrow D$, gives $\mu(B(m,r)) \geq \vol_K(r)/\vol_K(D)$. The ratio $\vol_K(r)/r^n$ is nonincreasing if $K \geq 0$ and nondecreasing if $K < 0$, so its infimum over $(0, D]$ is $\vol_K(D)/D^n$ or $v_n$, respectively.
\end{proof}

Combining \Cref{Prop:Ahlfors-Riemannian} with \Cref{Cor:Voronoi-Stability} gives the following generalization of \Cref{Thm:Main-Intro}, with \Cref{Thm:Main-Intro} being the $K=0$ case.

\begin{theorem}\label{Thm:Ricci}
    Let $M$, $\mu$, $\delta$, $X$, $Y$, $\omega_X$, $\omega_Y$, and $r$ be as in \Cref{Thm:Main-Intro}, except that we only assume $\Ric_M \geq (n-1)K$ for some $K \in \R$. Then
    \[
      d_\infty^{W_1}(\varphi_X^r, \varphi_Y^r)
      \le \Lambda\, \dTr(\omega_X, \omega_Y)
      \le 3r\Lambda\, \Wass(\omega_X, \omega_Y)
      \le 6r\delta\Lambda,
    \]
    where $\Lambda := (n+1)D^n/(r-\delta)^{n+1}$ if $K \geq 0$ and $\Lambda := (n+1)\vol_K(D)/(v_n (r-\delta)^{n+1})$ if $K < 0$.
\end{theorem}

\begin{proof}
    By \Cref{Prop:Ahlfors-Riemannian}, $\Lambda = (n+1)/(c(r-\delta)^{n+1})$ in both cases, so \Cref{Cor:Voronoi-Stability} applies with $\beta = 1/\Lambda$, giving the first inequality. The second follows from \Cref{Prop:Localized-Metric-Properties}\,(i). For the third, the coupling $(p_X, p_Y)_\# \mu$ moves each point a distance of at most $2\delta$, so $\Wass (\omega_X, \omega_Y) \leq 2 \delta$.
\end{proof}

We note that while Bishop--Gromov yields lower Ahlfors regularity constants for any compact Riemannian manifold, when the volume of small balls is known exactly \cite{GraySmallBalls}, applying \Cref{Prop:Ball-Mass-Stability} directly can yield significantly tighter estimates. To make the relevant constants completely explicit, below we take $r = \conv(M)/2$, i.e., half the \emph{convexity radius} of $M$ \cite[Section~6.4]{petersenRiemannianGeometry2016}, so that the $r$-balls form a good cover of $M$. In particular, \cite[Theorem~2.6]{dibbleConvexityRadius2017} gives $\conv (\mathbb R^n/\mathbb Z^n) = 1/4$ and $\conv (S^n)= \pi/2$.

\begin{example}\label{Ex:FlatNTorus}
On the flat torus $M = \mathbb{R}^n / \mathbb{Z}^n$, we have $r = 1/8$; assume $\delta < r$. Since $r$-balls are Euclidean for this value of $r$, we have $\mu(B(m,s)) = v_n s^n$ for $s \leq r$. Exactly as in the proof of \Cref{Cor:Voronoi-Stability}, we get $\omega_X(B(m,t)) \geq v_n(t-\delta)^n$ for $\delta<t\leq r$, and likewise for $\omega_Y$. Thus, \Cref{Prop:Ball-Mass-Stability} applies with
\[
\beta = \int_0^{1/8} v_n |t-\delta|_+^n \,dt = \frac{v_n}{n+1} \left( \frac 18- \delta\right)^{n+1}.
\]
Since $\dTr(\omega_X, \omega_Y) \leq 6r\delta = 3\delta/4$, we obtain
\[
d_\infty^{W_1}(\varphi_X^r, \varphi_Y^r)
\leq \frac{3(n+1)\delta}{4v_n \left( \frac 18 - \delta \right)^{n+1}}.
\]
The same argument on the unit sphere $S^n$, with $r = \pi/4$ and $\mu(B(m,s)) = A_n^{-1}\int_0^s \sin^{n-1}u\,du$, where $A_n := \int_0^\pi \sin^{n-1}u\,du$, gives
\[
d_\infty^{W_1}(\varphi_X^r, \varphi_Y^r)
\leq \frac{3\pi\delta A_n}{2\int_0^{\pi/4-\delta}\left(\frac{\pi}{4}-\delta-u\right)\sin^{n-1}u\,du}.
\]
\end{example}

\subsection{Convergence}
\label{Sec:Convergence}

To conclude this section, we show that if the weights converge in the $1$-Wasserstein distance to a lower Ahlfors regular reference measure, then \Cref{Thm:Main-Theorem} implies that their generalized nerve maps converge as well; in particular, this holds for Voronoi weights on a sequence of finite landmark sets converging to $M$ in the Hausdorff distance.

Let $\mu \in \mathcal P(M)$ be lower Ahlfors $(c,q)$-regular, and fix $r \in (0, D]$. Then every ball has positive $\mu$-mass, so $\supp \mu = M$, and by \Cref{Lem:Layer-Cake},
\begin{equation}\label{Eq:Reference-Mass}
S_\mu^r(m) = \int_0^r \mu(B(m,t))\,dt \geq \int_0^r c\,t^q\,dt = \frac{c}{q+1}\,r^{q+1} =: a_r \qquad \text{for every } m \in M.
\end{equation}
Thus the limiting nerve map $\varphi_\mu^r : M \to \mathcal P(M)$ is defined on all of $M$. Note that its values need not have finite support, so its image need not lie in $\check C^m (M; r)$. However, the inclusion $\supp\varphi_\mu^r(m)\subseteq\overline B(m,r)$ implies that its image lies in the \emph{closed $\infty$-\v Cech metric thickening}, defined by
\[
\left\{\alpha\in\mathcal P(M):\inf_{z\in M}\sup_{a\in\supp\alpha}d_M(z,a)\leq r\right\}.
\]
This is the closed analogue of \cite[Definition~3.5]{adams_persistent_2024}.

\begin{corollary}\label{Cor:Convergence}
    Let $\nu_k \in \mathcal P(M)$, $k\in \mathbb{N}$, satisfy $\Wass(\nu_k, \mu) \to 0$. Then for all $k$ with $\Wass(\nu_k, \mu) < a_r$, the map $\varphi_{\nu_k}^r$ is defined on all of $M$, and
    \[
    \sup_{m \in M} \Wass\big(\varphi_{\nu_k}^r(m), \varphi_\mu^r(m)\big) \leq \frac{3(q+1)}{c\,r^q}\, \Wass(\nu_k, \mu).
    \]
    Thus, $\varphi_{\nu_k}^r \to \varphi_\mu^r$ uniformly on $M$.
\end{corollary}

\begin{proof}
    Fix $k$ with $\Wass(\nu_k, \mu) < a_r$, fix $m \in M$, and let $\pi \in \Pi(\nu_k, \mu)$. Since $k_r(m, \cdot)$ is $1$-Lipschitz,
    \[
    \big| S_{\nu_k}^r(m) - S_\mu^r(m) \big| \leq \int \big| k_r(m,x) - k_r(m,z) \big|\, d\pi(x,z) \leq \int d_M \, d\pi,
    \]
    and infimizing over $\pi$ gives $S_{\nu_k}^r(m) \geq a_r - \Wass(\nu_k,\mu) > 0$, so $\varphi_{\nu_k}^r$ is defined at $m$. By \eqref{Eq:Reference-Mass}, $\max\{S_{\nu_k}^r(m), S_\mu^r(m)\} \geq a_r$ for every $m$, so $\beta_\ast(\nu_k, \mu) \geq a_r$. \Cref{Thm:Main-Theorem} with $\beta = a_r$ then gives the bound, since $3r/a_r = 3(q+1)/(c\,r^q)$.
\end{proof}

For example, let $X_k \subseteq M$ be finite sets with 
\[
\delta_k := \sup\limits_{z \in M} d_M(z, X_k) \to 0,
\]
let $p_k : M \to X_k$ be nearest neighbor projections, and set $\omega_k := (p_k)_\#\mu$. Every landmark has a Voronoi neighborhood of positive $\mu$-measure, so $\supp \omega_k = X_k$, and the coupling $(p_k, \id)_\#\mu$ gives $\Wass(\omega_k, \mu) \leq \delta_k$. Hence \Cref{Cor:Convergence} gives, for all $k$ large enough,
\[
d_\infty^{W_1}(\varphi_{X_k}^r, \varphi_\mu^r) \leq \frac{3(q+1)}{c\,r^q}\,\delta_k,
\]
and $\varphi_{X_k}^r \to \varphi_\mu^r$ uniformly on $M$.

\section{Experiments}
\label{Sec:Experiments}

In this section, we test the sharpness and computational cost of the three estimates of \Cref{Thm:Main-Theorem} through numerical experiments. First, we consider six real and synthetic biological datasets, using them to illustrate different methods of equipping an unweighted point cloud with landmark weights, and examine how different choices affect the sharpness of our estimates (\Cref{Sec:Bio-Data}). We then test how the sharpness of the localized transport estimate varies with the dimension for uniform samples of two families of manifolds (\Cref{Sec:Manifold-Experiments}). Finally, we compare the time and memory required to compute each estimate exactly on samples of the unit square of increasing size (\Cref{Sec:Computational-Cost}). 

For brevity, we write throughout this and subsequent sections
\[
Q_1 := d_\infty^{W_1} (\varphi_X^r, \varphi_Y^r), \quad Q_{\mathrm{pt}} := \Delta_L^r(\omega_X, \omega_Y), \quad Q_{\mathrm{tr}} := \frac{d_L^r(\omega_X, \omega_Y)}{\beta_\ast}, \quad Q_{\mathrm{was}} := \frac{3r}{\beta_\ast} \Wass(\omega_X, \omega_Y).
\]
Recall that we refer to $Q_{\mathrm{pt}}$ as the \emph{pointwise estimate}, and to $Q_{\mathrm{tr}}$ and $Q_{\mathrm{was}}$ as the \emph{localized transport} and \emph{Wasserstein} estimates, respectively. The chain of inequalities of \Cref{Thm:Main-Theorem} is then 
\[
Q_1 \leq Q_{\mathrm{pt}} \leq Q_{\mathrm{tr}} \leq Q_{\mathrm{was}},
\]
and a ratio $Q_\bullet/Q_1$ closer to $1$ indicates a sharper estimate.

\subsection{Biological Data}
\label{Sec:Bio-Data}

Here we aim to investigate how choices of landmark weights affect the sharpness of \Cref{Thm:Main-Theorem}. We use six real and synthetic datasets from Scanpy \cite{scanpyDatasets}, a Python toolkit for analyzing single-cell gene expression data, and compare three weighting schemes: uniform, Voronoi, and weights obtained by maximizing the baseline lower bound
\[
\min\limits_{m\in M}S_X^r(m)\leq \beta_\ast.
\]

\begin{table}[htbp]
\centering\small
\caption{Datasets and landmark selection. Distances are Euclidean in the indicated representation, and $D=\diam(M)$. Here $\#L$ is the max-min candidate count and $\#X_{\rm opt}$ is the optimized support size. Uniform and Voronoi weights use all of $L$.}
\label{tab:weighting-setup}
\setlength{\tabcolsep}{5pt}
\begin{tabular}{@{}lrlrrrr@{}}
\toprule
Dataset & $\#M$ & Representation & $\#L$ & $\# X_{\rm opt}$ & $r$ & $r/D$\\
\midrule
\texttt{toggleswitch} & 200 & Original & 20 & 20 & 0.19664 & 0.1423\\
\texttt{krumsiek11} & 640 & Original & 64 & 64 & 0.20724 & 0.0783\\
\texttt{moignard15} & 3934 & Standardized & 100 & 95 & 8.266 & 0.2492\\
\texttt{paul15} & 2730 & Preprocessed PCA & 100 & 99 & 27.722 & 0.2127\\
\texttt{pbmc68k\_reduced} & 700 & Existing PCA & 100 & 97 & 13.243 & 0.2716\\
\texttt{pbmc3k\_processed} & 2638 & Existing PCA & 100 & 95 & 18.480 & 0.2200\\
\bottomrule
\end{tabular}
\end{table}

For \texttt{moignard15}, we standardize each coordinate to have zero mean and unit variance. For \texttt{paul15}, each data point represents a cell whose coordinates are gene-expression counts. We preprocess the data by rescaling each point to have $\ell^1$ norm $10^4$, applying $x\mapsto \log (1+x)$ coordinatewise, standardizing each coordinate to have zero mean and unit variance across the dataset, and projecting onto the first $30$ principal coordinates. For the two PBMC datasets, we use the first $30$ supplied principal components. We then do a deterministic max-min run, initialized with seed $0$, selecting the candidates $L$ once per dataset. The scale $r$ is then chosen  for each dataset after inspecting coverage and neighborhood size, as well as exploratory persistence. The radii were about $2\Cov (X)$ for the first two datasets and about $1.2 \Cov(X)$ for the others.

\begin{figure}[!tp]
\centering
\begin{minipage}[c]{.27\linewidth}
  \centering
  \includegraphics[width=\linewidth]{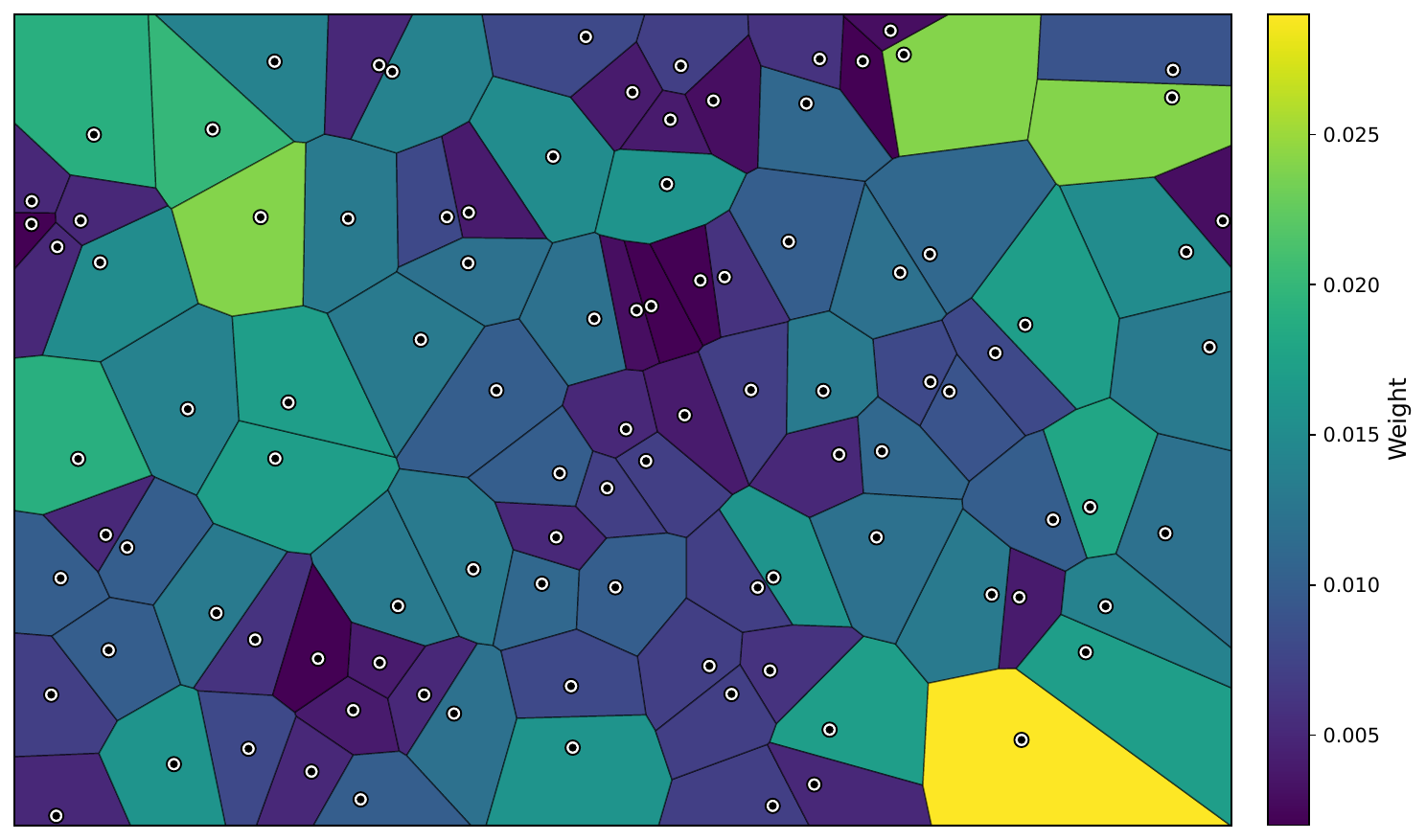}
\end{minipage}\hfill
\begin{minipage}[c]{.71\linewidth}
  \centering
  \includegraphics[width=\linewidth]{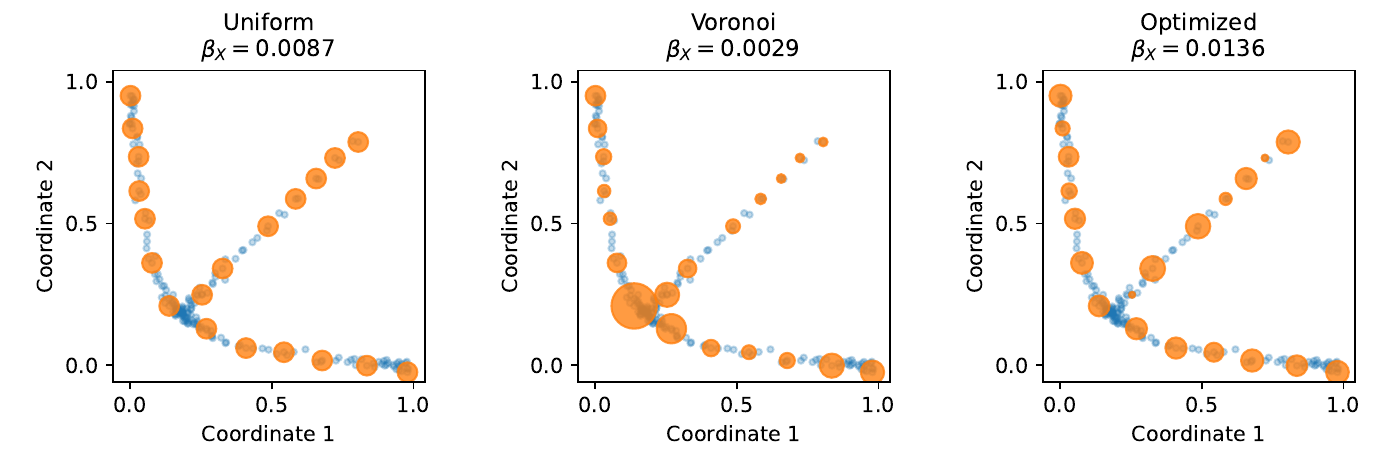}
\end{minipage}
\caption{Left: An illustration of a Voronoi tessellation with cells colored according to their measure. 
Right: landmarks and weights on \texttt{toggleswitch}, under uniform, Voronoi, and optimized weights. Data points are blue and landmark disk areas are proportional to their weights (normalized to be on a common scale).}
\label{Fig:Voronoi-Weights}
\label{fig:weighting-landmarks}
\end{figure}

The three weighting schemes we use are as follows. With 
\[
\mu_M := \frac{1}{\card M} \sum\limits_{m \in M} \delta_m,
\]
the first two schemes (uniform and Voronoi) take $X=L$ and weights to be 
\[
\frac{1}{\card L} \sum\limits_{x\in L} \delta_x \quad \text{and} \quad (p_L)_\# \mu_M,
\]
respectively. In other words, Voronoi weights count the number of observations assigned to each landmark. For the third scheme, we use a linear program to maximize
\[
b_L := \max_{ \omega_x \geq 0,\; \sum \omega_x=1} \min_{m \in M} \sum_{x \in L} \omega_x k_r(m, x).
\]
Numerically, we minimize support size and then covering radius subject to a baseline normalizer of at least $b_L-10^{-8}$. After reoptimizing on this support, we mix the weights with the uniform measure in proportion $10^{-6}$ to guarantee the weights are all positive. We call the resulting weights \emph{$\beta$-optimal weights.} \Cref{fig:weighting-landmarks} illustrates the three weight choices on the same max-min landmark run.

We use the following perturbation model in this experiment. In each trial, we draw independent standard Gaussian vectors $Z_i$, and set
\[
y_i := p_M (x_i + \epsilon \Cov(X) Z_i ), \quad f(x_i) := y_i, \quad \omega_Y := f_\# \omega_X,
\]
where $p_M$ denotes the nearest neighbor in $M$. To ensure that perturbed points don't return to their original position after projecting back onto $M$, we separately calibrate $\epsilon$ so that approximately $5\%, 10\%, 20\%$, or $40\%$ of the initial landmarks move. We solve all transport problems using POT's network simplex implementation \cite{flamary2021pot}, aggregating mass outside each $r$-ball onto a single point to streamline the computations (see \Cref{Sec:Computational-Cost}).  The results are displayed in \Cref{tab:weighting-results}, where each entry represents the mean over ten trials $\pm$ one standard deviation.

\begin{table}[!htbp]
\centering\small
\caption{Weighting experiment at the $40\%$ movement target: means $\pm$ sample standard deviations over ten trials. $D=\diam(M)$; the observed distance is reported as a percentage of $D$.}
\label{tab:weighting-results}
\setlength{\tabcolsep}{4pt}
\begin{tabular}{@{}llrrrr@{}}
\toprule
Dataset & Weights & $100Q_1/D$ & $Q_{\rm pt}/Q_1$ & $Q_{\rm tr}/Q_1$ & $Q_{\rm was}/Q_1$\\
\midrule
\texttt{toggleswitch} & Uniform & $2.61\pm0.35$ & $2.25\pm0.25$ & $4.7\pm1.1$ & $15.3\pm2.5$ \\
& Voronoi & $3.26\pm0.65$ & $2.29\pm0.22$ & $26.8\pm10.6$ & $48.8\pm14.7$ \\
& Optimized & $2.58\pm0.58$ & $2.23\pm0.14$ & $3.1\pm0.3$ & $10.2\pm1.7$ \\
\addlinespace[2pt]
\texttt{krumsiek11} & Uniform & $2.01\pm0.33$ & $2.34\pm0.23$ & $5.0\pm0.9$ & $43.6\pm7.1$ \\
& Voronoi & $2.31\pm0.62$ & $2.26\pm0.30$ & $20.7\pm3.8$ & $109\pm21$ \\
& Optimized & $2.13\pm0.29$ & $2.32\pm0.21$ & $3.0\pm0.4$ & $26.1\pm3.1$ \\
\addlinespace[2pt]
\texttt{moignard15} & Uniform & $22.48\pm2.54$ & $1.67\pm0.14$ & $29.5\pm6.5$ & $421\pm52$ \\
& Voronoi & $24.43\pm1.88$ & $1.70\pm0.11$ & $1645\pm447$ & $11136\pm3030$ \\
& Optimized & $24.01\pm1.84$ & $1.72\pm0.09$ & $10.5\pm0.9$ & $125\pm8$ \\
\addlinespace[2pt]
\texttt{paul15} & Uniform & $20.65\pm0.89$ & $1.86\pm0.14$ & $38.5\pm9.3$ & $422\pm72$ \\
& Voronoi & $22.36\pm1.62$ & $1.89\pm0.28$ & $1138\pm270$ & $6919\pm1226$ \\
& Optimized & $21.40\pm1.63$ & $1.76\pm0.22$ & $14.6\pm2.0$ & $146\pm11$ \\
\addlinespace[2pt]
\texttt{pbmc68k\_reduced} & Uniform & $27.20\pm1.93$ & $1.94\pm0.14$ & $34.4\pm11.2$ & $455\pm116$ \\
& Voronoi & $28.81\pm3.71$ & $1.88\pm0.29$ & $298\pm46$ & $1890\pm343$ \\
& Optimized & $28.76\pm2.59$ & $1.74\pm0.17$ & $14.4\pm2.3$ & $159\pm19$ \\
\addlinespace[2pt]
\texttt{pbmc3k\_processed} & Uniform & $20.85\pm1.70$ & $1.82\pm0.16$ & $40.0\pm10.3$ & $335\pm54$ \\
& Voronoi & $22.29\pm2.78$ & $1.69\pm0.22$ & $1244\pm560$ & $4909\pm1486$ \\
& Optimized & $21.14\pm2.14$ & $1.84\pm0.18$ & $14.8\pm2.9$ & $127\pm6$ \\
\bottomrule
\end{tabular}
\end{table}

Across all $72$ dataset weighting and movement combinations, the average $Q_{\mathrm{pt}} / Q_1$ ranges from $1.34$ to $2.34$. The other estimates have comparatively wide ranges: the mean $Q_{\mathrm{tr}}/Q_1$ ranges from $2.6$ to $14.8$ for optimized weights, $4.0$ to $40.0$ for uniform weights, and about $13.6$ to $1681$ for Voronoi weights. The optimized weights exhibited the best quantitative stability among the three weighting choices, while Voronoi weights fared the least well. This is partially explained by the fact that Voronoi weights appear to sometimes assign microscopic mass to sparse regions, making $\beta_\ast$ extremely small. The choices of weights thus offer tradeoffs between faithfulness to the underlying distribution and quantitative stability.

\subsection{Manifold Data}
\label{Sec:Manifold-Experiments}

To study how the sharpness of the localized transport estimate varies with the dimension, we use samples from two families of manifolds. For each $n=2,\ldots,10$, we sample $250\cdot2^n$ points uniformly from the flat torus $(\R/2\pi\Z)^n$ and the unit sphere $S^n$, using geodesic distances. A max-min run selects $X$, and a second run restricted to $M\setminus X$ selects $Y$, both with $\card M / 10$ point, and we set
\[r=\min\{2\max\{\Cov(X),\Cov(Y)\},D\}.\]
Since the samples are uniform and max-min landmarks are well separated, the Voronoi weights of \Cref{Thm:Ricci} are approximately uniform, so for simplicity we use uniform weights. We also evaluate $\beta_\ast$ directly on $M$, rather than using the lower bound of \Cref{Thm:Ricci}, since the Bishop--Gromov lower bound on ball volumes underlying it is far from tight, as \Cref{Ex:FlatNTorus} illustrates.

\begin{figure}[!tp]
\centering
\includegraphics[width=.78\linewidth]{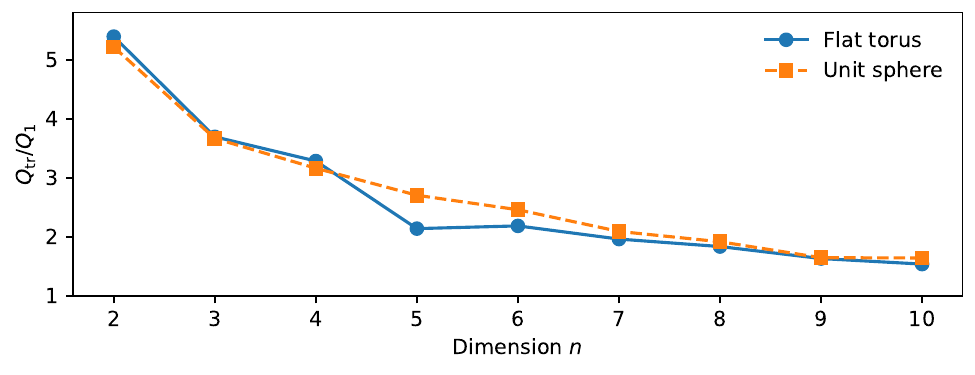}
\caption{Sharpness of the localized transport estimate on sampled manifolds. Each curve uses uniform landmark weights and the empirically evaluated $\beta_*$.}
\label{Fig:Empirical-Comparison}
\end{figure}

The ratio $Q_{\rm tr}/Q_1$ decreases overall from about $5.4$ and $5.2$ at $n=2$ to $1.54$ and $1.65$ at $n=10$ for the torus and sphere, respectively (\Cref{Fig:Empirical-Comparison}). Empirically, both ratios appear to tend toward $1$ as the dimension grows.

\subsection{Computational Cost of the Estimates}
\label{Sec:Computational-Cost}

We briefly address the computational efficiency of calculating each of the three estimates of \Cref{Thm:Main-Theorem}. For the first two estimates (i.e., the local estimates), since each localized cost function vanishes outside $B(m, r)$, all mass outside $B(m,r)$ can be merged into a single auxiliary point (one for each of $X$ and $Y$) without changing the transport problem. Thus computing both local estimates requires solving $|M|$ transportation problems with at most $(k+1)^2$ variables each, where
\[
 k=\max_{m\in M}\max\{\card (X\cap B(m,r)), \card (Y\cap B(m,r))\}.
\]
Given the neighborhoods $X\cap B(m, r)$ and $Y \cap B(m,r)$\footnote{Such neighborhoods can be found efficiently using, e.g., cover trees \cite{beygelzimer2006cover}.}, the time and storage required to compute both estimates are $O(\card M)$ for bounded $k$; this holds, e.g., for max-min landmarks in metric spaces of bounded doubling dimension \cite{gupta2003bounded} when the scale is comparable to the covering radii, by a packing argument \cite{smid2009weak}. On the other hand, the exact transport solvers used here to calculate $\Wass$ instead require the full $(\card X)\times(\card Y)$ cost matrix. Entropic regularization \cite{cuturi2013sinkhorn} offers a faster approximation of Wasserstein distances, and computes an $\epsilon$-accurate transport cost in time nearly linear in the number of cost-matrix entries \cite{altschuler2017near}. However, such regularization still operates on the full $(\card X) \times (\card Y)$ cost matrix, whereas each local problem involves only at most $(k+1)^2$ entries.

To test this empirically, we sample $N=2000\cdot2^j$ points uniformly from $[0,1]^2$, $j=0,\ldots,7$, select two random-start max-min sets of size $N/10$, and set $r=2\max\{\Cov(X),\Cov(Y)\}$. Tests use uniform or empirical Voronoi weights, POT's network simplex solver, and, for uniform weights, SciPy's assignment solver. Each estimate runs in a separate single-core process on a computer with 32 GB of memory. Note also that the reported times do not include landmark selection.

\begin{figure}[!tp]
\centering
\includegraphics[width=.85\linewidth]{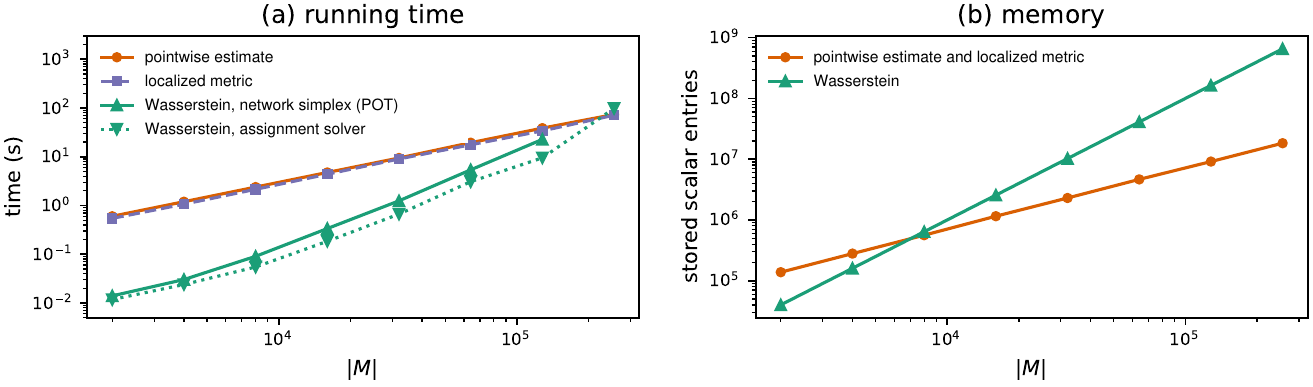}
\caption{Exact-computation timings for uniform weights (left) and stored scalar entries (right). The network simplex curve stopped at $N=128{,}000$ due to insufficient memory.}
\label{fig:computational-cost}
\end{figure}

For $k\leq12$, computing both local estimates scales approximately linearly. The assignment solver is faster through $N=128{,}000$ (23 versus 39 seconds for the pointwise estimate), but at $N=256{,}000$ takes 98 seconds and 5.0 GB, versus 71--74 seconds and 0.22 GB for the local estimates. The network simplex run fails at that size for both weightings due to insufficient memory. Thus, the two sharper estimates of \Cref{Thm:Main-Theorem} are less costly to compute on large datasets than the exact Wasserstein estimate; see \Cref{fig:computational-cost}.

\section{Applications to Topological Coordinates}
\label{Sec:Applications}

Many coordinate pipelines take the form $\Psi_X := F_X \circ \varphi_X^r$ with $F_X : \check C^m(X; r) \to Z$ for a metric space $(Z, d_Z)$. To describe the role of nerve map stability in the stability of such pipelines, we observe for example that if $F_X$ and $F_Y$ are restrictions of a common $L$-Lipschitz map $F : A \to Z$, where $A \subseteq \check C^m(M; r)$ contains the images of $\varphi_X^r$ and $\varphi_Y^r$, then
\begin{equation}\label{Eq:Application-Stability}
\sup_{m \in X^r \cap Y^r} d_Z\big(\Psi_X(m), \Psi_Y(m)\big) \leq L \, d_\infty^{W_1}(\varphi_X^r, \varphi_Y^r),
\end{equation}
which \Cref{Thm:Main-Theorem} bounds in turn. One way to obtain such maps is via affine extension: for a normed vector space $V$, $A \subseteq M$, and $F : A \to V$, set 
\[
\hat F(\alpha) := \sum\limits_{a \in A} \alpha(a) F(a)
\]
for finitely supported $\alpha \in \mathcal P(M)$ with $\supp \alpha \subseteq A$.

\begin{lemma}\label{Lem:Lipschitz-Composite}
    If $F : A \to V$ is $L$-Lipschitz, then $\| \hat F(\alpha) - \hat F(\alpha') \|_V \leq L \, \Wass(\alpha, \alpha')$ for all such $\alpha, \alpha'$.
\end{lemma}

\begin{proof}
    For $\pi \in \Pi(\alpha, \alpha')$, we have $\| \hat F(\alpha) - \hat F(\alpha') \|_V = \| \int (F(a) - F(b)) \, d\pi(a,b) \|_V \leq L \int d_M \, d\pi$; now infimize over $\pi$.
\end{proof}

\subsection{Out-of-Sample Extension and Classification}
\label{Sec:Classification-Experiment}

Given landmark coordinates $F:X\to V$, the composite $\Psi_X:=\hat F\circ\varphi_X^r$ assigns out-of-sample coordinates on $X^r$ by averaging the landmark values according to the nerve map weights; see \cite{bengio2004out,williams2001nystrom} for background on out-of-sample methods. Note that this procedure need not preserve $F$ on the landmarks. In the special case when $F_X$ and $F_Y$ are restrictions of a common $L$-Lipschitz $F_M : M \to V$, \Cref{Lem:Lipschitz-Composite} gives \Cref{Eq:Application-Stability} with $Z = V$. 

The goal of this experiment is to illustrate how stability of the nerve map controls the sensitivity of out-of-sample classification scores under landmark perturbation, and to examine how such perturbations affect the predicted labels. For the classification task, we use the Wisconsin Diagnostic Breast Cancer dataset \cite{wolberg1993breast,scikit-learn}, with a stratified 70/30 split (random seed $17$) into $398$ training and $171$ test observations. We standardize each feature using training-set means and standard deviations, and fit the two-component PCA on the standardized training observations. The fitted transformations are applied to the test observations. The two components retain $64.29\%$ of the training variance. A linear SVM in this plane supplies the $1$-Lipschitz score $F(z)=(w^\top z+b)/\|w\|_2$. We select eight max-min training landmarks with uniform weights, which are transferred to the perturbed centers by the landmark correspondence. Let $M_{\rm train}$ denote the training observations in the PCA plane, and set
\[
c_X:=\max_{z\in M_{\rm train}}\min_{x\in X}\|z-x\|_2,
\qquad r=2c_X=10.2098.
\]
Define
\[
 \hat f_X(m)=\sum_{x\in X}\phi_x^r(m)F(x),\qquad
 \hat y_X(m)=\mathbf1_{\{\hat f_X(m)\geq0\}}.
\]
In each trial, landmarks move by $\epsilon c_X$ in independent random unit directions, with the score function and radius fixed. The background space is the union of test observations and both landmark sets. By \Cref*{Lem:Lipschitz-Composite}, the score error satisfies 
\[
Q_S:=\max\limits_{X^r\cap Y^r}|\hat f_X-\hat f_Y|\leq Q_1\leq Q_{\rm pt}.
\]

\begin{table}[!htbp]
\centering\small
\caption{Classification under landmark perturbation. Flips are relative to baseline predictions; baseline test accuracy is $98.25\%$.}
\label{tab:classification}
\setlength{\tabcolsep}{4pt}
\begin{tabular}{@{}rrrrrr@{}}
\toprule
$100\epsilon$ & $Q_S$ & $Q_1$ & $Q_{\rm pt}/Q_1$ & Flips (\%) & Accuracy (\%)\\
\midrule
2 & $0.113\pm0.023$ & $0.234\pm0.014$ & $2.38\pm0.17$ & $0.82\pm0.56$ & $98.36\pm0.37$\\
5 & $0.279\pm0.057$ & $0.580\pm0.036$ & $2.39\pm0.20$ & $1.52\pm1.07$ & $97.89\pm1.07$\\
10 & $0.554\pm0.112$ & $1.123\pm0.075$ & $2.47\pm0.22$ & $2.34\pm1.72$ & $96.96\pm1.58$\\
20 & $1.077\pm0.205$ & $2.094\pm0.117$ & $2.55\pm0.21$ & $3.63\pm1.80$ & $95.67\pm1.70$\\
\bottomrule
\end{tabular}
\end{table}

The score error and nerve-map distance grow approximately linearly with the perturbation, while $Q_{\rm pt}/Q_1$ remains between $2.4$ and $2.6$. At the largest perturbation, mean accuracy is $95.67\%$ and $3.63\%$ of labels change (\Cref{tab:classification}).

\subsection{Laplacian Eigenvector Cascading}
\label{Sec:Cascading-Experiment}

The Laplacian eigenvector cascading algorithm of \cite{MikePerea2019LaplacianEigenvectorCascading} computes graph Laplacian eigenvectors on a cover tower of increasingly fine landmark sets, and extends these eigenvectors to the entire dataset using nerve maps. This procedure has two types of instability: the nerve map changes between levels, which a two-scale version of \Cref{Thm:Main-Theorem} and \Cref{Lem:Lipschitz-Composite} controls, as we show below, and also the eigenvectors selected at consecutive scales may disagree, which the cascade algorithms are designed to control. Here we isolate the former source of instability by holding the eigenvector block fixed and comparing its coordinate extensions through the nerve maps at consecutive scales.

For this experiment, we use the PIN\_1 and SPHERE examples of \cite[Examples~4 and~5]{MikePerea2019LaplacianEigenvectorCascading}, consisting of $15{,}148$ points in three triangles and $12{,}432$ points in six spherical annuli, with Euclidean and geodesic distances, respectively. In each of ten trials, a random-start max-min ordering gives nested sets $C_0\subset C_1\subset C_2\subset C_3$ with $\Cov(C_i)\leq\epsilon_0 2^{-i}$, where $\epsilon_0=0.15$ for PIN\_1 and $0.35$ for SPHERE. Unlike the original construction, we use nearest-neighbor maps $\tau_i:C_{i+1}\to C_i$ fixing $C_i$. With $\delta_i=\max_y d_M(y,\tau_i(y))$, set $r_3=2\Cov(C_3)$ and $r_i=1.05(r_{i+1}+\delta_i)$. The landmark weights are $\omega_{C_3}=(p_{C_3})_\#\mu_M$ and $\omega_{C_i}=(\tau_i)_\#\omega_{C_{i+1}}$, where $\mu_M$ is uniform on $M$.

The edge weights we use here are $W_i(a,b)=\omega_{C_i}(a)\omega_{C_i}(b)k_{2r_i}(a,b)$ for $a\neq b$, with zero diagonal. We construct spectral coordinate blocks $F_i:C_i\to\R^k$ using the second cascade of \cite{MikePerea2019LaplacianEigenvectorCascading}: at the coarsest level we use eigenvectors corresponding to the $k$ smallest positive eigenvalues of $I-D_i^{-1}W_i$, where $D_i$ is the degree matrix; at finer levels we align the pulled-back coordinates within the corresponding eigenvalue-cluster spans, using relative clustering threshold $0.25$. This groups neighboring positive eigenvalues when the larger is less than $25\%$ above the smaller; differences below $10^{-9}$ are also treated as ties. We normalize the blocks by $F_i^\top D_iF_i=(\operatorname{tr}D_i)I_k$, with $k=2$ for PIN\_1 and $3$ for SPHERE. Common orthogonal changes of the output basis do not affect the errors below.

For consecutive levels, put $X=C_{i+1}$, $Y=C_i$, $r=r_{i+1}\leq s=r_i$, and $\tau=\tau_i$. Holding $F_i$ fixed, compare
\[
f_Y=\widehat F_i\circ\varphi_Y^s,\qquad
f_X^\tau=\widehat F_i\circ\tau_\#\circ\varphi_X^r,\qquad
E_\tau=\max_{m\in M}\|f_Y(m)-f_X^\tau(m)\|_2.
\]
Let $A_m=\supp\varphi_X^r(m)\cup\supp\varphi_Y^s(m)$, $a_m=\max\{S_X^r(m),S_Y^s(m)\}$, and $\beta_*^{r,s}=\min_m a_m$. The two-scale costs are
\[
c_m^{r,s}(x,y)=\min\{k_r(m,x),k_s(m,y)\}\,d_M(x,y)
 +(r+s)|k_r(m,x)-k_s(m,y)|,
\]
with $\widetilde c_m^{\,r,s}$ obtained by replacing $r+s$ with $\min\{\diam A_m,r+s\}$. Put
\[
q_1(m)=\Wass(\varphi_X^r(m),\varphi_Y^s(m)),\quad
q_{\rm pt}(m)=\frac{T_{\widetilde c_m^{\,r,s}}(\omega_X,\omega_Y)}{a_m},\quad
q_{\rm tr}(m)=\frac{T_{c_m^{r,s}}(\omega_X,\omega_Y)}{\beta_*^{r,s}}.
\]
Writing $L(m)=\operatorname{Lip}((F_i\circ\tau)|_{A_m})$ and $L_*=\max_m L(m)$, let $B_j=\|Lq_j\|_\infty$ for $j=1,\mathrm{pt},\mathrm{tr}$, and
\[
B_{\rm was}=\frac{L_*}{\beta_*^{r,s}}
\big[(2r+s)\Wass(\omega_X,\omega_Y)+(r+s)(s-r)\big].
\]
Indeed, $c_m^{r,s}\leq(2r+s)d_M+(r+s)(s-r)$, and the proof of \Cref{Thm:Main-Theorem}, applied at two scales, together with \Cref{Lem:Lipschitz-Composite}, gives
\[
E_\tau\leq B_1\leq B_{\rm pt}\leq B_{\rm tr}\leq B_{\rm was}.
\]

\begin{table}[!htbp]
\centering\small
\caption{Cascading errors and bound ratios. $E_{\rm rel}$ and $\overline E_{\rm rel}$ are the maximum and RMS errors divided by $\diam f_Y(M)$, expressed as percentages. Transitions are between levels $i$ and $i+1$. Entries summarize ten trials.}
\label{tab:cascade-ratios}
\setlength{\tabcolsep}{3pt}
\begin{tabular}{@{}lcrrrrrr@{}}
\toprule
Dataset & $i$ & $E_{\rm rel}$ & $\overline E_{\rm rel}$ & $B_1/E_\tau$ & $B_{\rm pt}/E_\tau$ & $B_{\rm tr}/E_\tau$ & $B_{\rm was}/E_\tau$\\
\midrule
PIN\_1 & 0 & $9.12\pm1.01$ & $1.83\pm0.21$ & $4.44\pm0.61$ & $14.8\pm2.5$ & $64.2\pm17.4$ & $966\pm276$\\
 & 1 & $8.74\pm2.19$ & $0.89\pm0.08$ & $4.06\pm0.96$ & $12.6\pm2.8$ & $64.0\pm22.4$ & $3700\pm980$\\
 & 2 & $6.57\pm1.07$ & $0.42\pm0.03$ & $4.72\pm1.09$ & $13.7\pm2.4$ & $53.9\pm13.6$ & $13900\pm3200$\\
\addlinespace[2pt]
SPHERE & 0 & $11.8\pm1.1$ & $4.84\pm0.23$ & $4.61\pm0.44$ & $15.2\pm1.7$ & $26.9\pm2.1$ & $436\pm39$\\
 & 1 & $8.18\pm0.88$ & $2.10\pm0.08$ & $4.23\pm0.65$ & $14.2\pm2.5$ & $31.2\pm5.2$ & $1820\pm300$\\
 & 2 & $5.67\pm0.87$ & $1.06\pm0.03$ & $3.77\pm0.60$ & $11.2\pm1.9$ & $29.6\pm2.8$ & $7010\pm850$\\
\bottomrule
\end{tabular}
\end{table}

\begin{figure}[!tp]
\centering
\includegraphics[width=.88\linewidth]{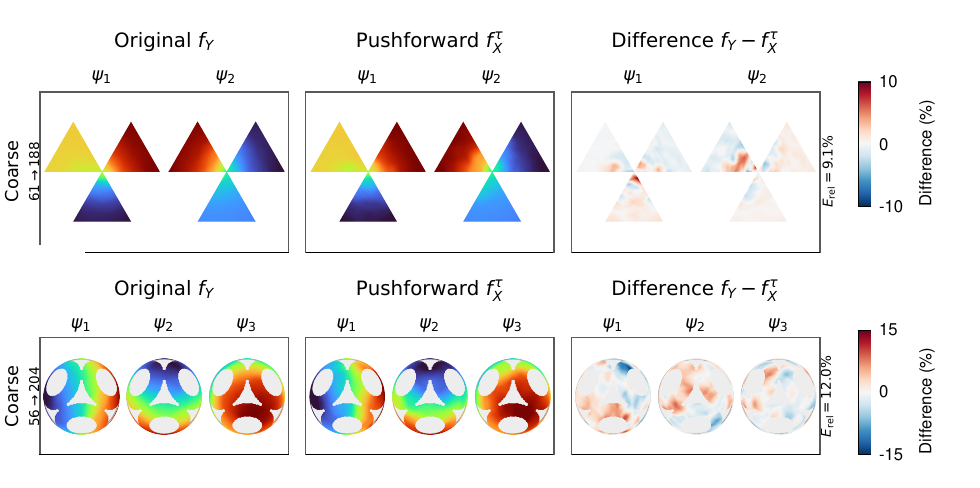}
\caption{Coarse-to-medium transitions for PIN\_1 and SPHERE. Each row shows the coarse coordinates, their pushforward extension, and componentwise differences as percentages of the coordinate range. All levels are summarized in \Cref{tab:cascade-ratios}.}
\label{fig:cascade-pin}
\label{fig:cascade-sphere}
\end{figure}

The normalized RMS error is below $5\%$ and approximately halves with each refinement. The pointwise estimate is about $11$--$15$ times the coordinate error, whereas the Wasserstein estimate exceeds it by factors of roughly $4\times10^2$ to $1.4\times10^4$ (\Cref{tab:cascade-ratios}). These comparisons hold the spectral data fixed and do not measure changes in independently recomputed eigenvectors.

\subsection{Eilenberg--MacLane Coordinates}
\label{Sec:SCC-Experiment}

Eilenberg--MacLane coordinates \cite{desilvaPersistentCohomology2011,pereaSparseCircular2020,pereaDREiMacDimensionality2023a,Perea2018MultiscaleProjective} take the form $\Psi_X = F_X \circ \varphi_X^r$, where $F_X : \check C^m(X; r) \to K(G,n)$ represents a persistent cohomology class of the \v{C}ech filtration on $X$, so \Cref{Eq:Application-Stability} applies whenever $F_X$ and $F_Y$ restrict a common Lipschitz map. In practice, however, they are built separately from cocycles on different complexes, so this need not hold. 

We aim to test the contribution of the nerve map in the stability of sparse circular coordinates, i.e., $K(\mathbb Z, 1)$-valued coordinates, isolating it by holding the cocycle and smoothing data fixed. We use $240$ noisy $32\times32$ images of a rotating arrow, blurred by one pixel and corrupted with Gaussian noise of standard deviation $0.01$. After clipping to $[0,1]$, centering, and normalizing each image, we use Euclidean distances in the first eight principal components, which retain $92.1\%$ of the variance. Take $30$ uniformly weighted landmarks $x_i=m_{8i}$ and $r=0.4$. For $q\in\{1,2,3,4\}$, let $y_i^{(q)}=m_{8i+q(-1)^i}$, with indices modulo $240$. We retain the labels and their weights $1/30$, adding masses at coincident centers; thus $Y_4$ has $15$ distinct centers.

Following the sparse circular coordinates construction of \cite{pereaSparseCircular2020}, we select the most persistent Vietoris--Rips $H^1$ class over $\mathbb F_{47}$ using Ripser \cite{Bauer2021Ripser}. Its interval is $[0.2837,0.8916)$, with $\max\{0.2837,\Cov(X)\}<r<0.8916/2$. The symmetric lift $\eta$ was verified to be an integer cocycle on $K=\operatorname{VR}_{2r}(X)$, whose simplices have pairwise distances less than $2r$. Here $d_{2r}:C^0(K;\mathbb R)\to C^1(K;\mathbb R)$ is the coboundary in the oriented simplex bases, and $+$ denotes its Moore--Penrose inverse for the standard unweighted Euclidean inner products. We regard $\eta$ as a real cochain in the smoothing formulas $\tau=-d_{2r}^+\eta$ and $\theta=\eta+d_{2r}\tau$. We use the circular-coordinate formula of \cite{pereaSparseCircular2020}, keeping $\tau$ and $\theta$ fixed on the original landmark labels while changing the partition of unity. Every perturbed support, interpreted on the original labels, spans a simplex of $K$, since its largest baseline edge length is $0.7163<2r$.

Let $h_X, h_{Y_q}: M \to S^1$ denote the resulting circular coordinates, where $Y_q := \{ y_i^{(q)} \}$. After aligning $h_X$ and $h_{Y_q}$ by a global phase minimizing their RMS geodesic distance, we record the maximum and RMS circular errors in \Cref{tab:scc-stability}, together with the nerve-map distances.

\begin{table}[!htbp]
\centering\small
\caption{Circular coordinates with cocycle and smoothing fixed. Landmark movement is measured in the PCA metric; circular errors are in radians.}
\label{tab:scc-stability}
\setlength{\tabcolsep}{7pt}
\begin{tabular}{@{}crrrrrr@{}}
\toprule
$q$ & $\max_i\|x_i-y_i^{(q)}\|_2$ & $Q_1$ & $Q_{\rm pt}$ & $Q_{\rm pt}/Q_1$ & $d_\infty^{S^1}$ & $d_{\rm RMS}^{S^1}$\\
\midrule
1 & 0.0420 & 0.0656 & 0.1873 & 2.86 & 0.0259 & 0.0197\\
2 & 0.0809 & 0.1216 & 0.3661 & 3.01 & 0.0537 & 0.0382\\
3 & 0.1111 & 0.1653 & 0.5119 & 3.10 & 0.0827 & 0.0547\\
4 & 0.1455 & 0.1953 & 0.6148 & 3.15 & 0.1073 & 0.0654\\
\bottomrule
\end{tabular}
\end{table}

\begin{figure}[!tp]
\centering
\includegraphics[width=.86\linewidth]{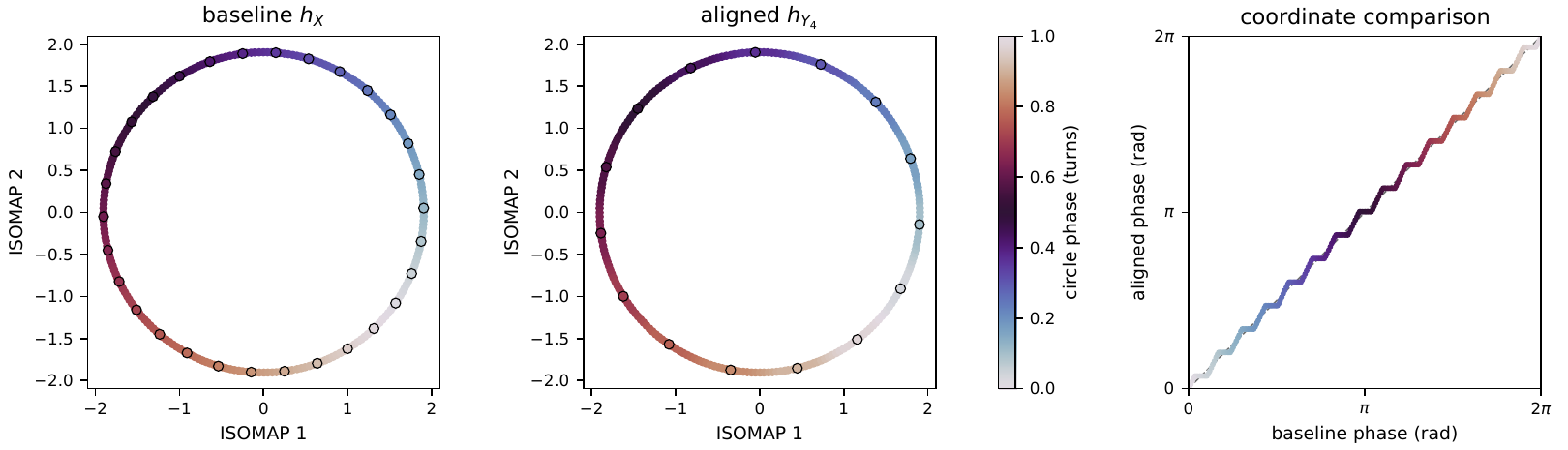}
\caption{Baseline and phase-aligned perturbed circular coordinates on an ISOMAP visualization, followed by their phase comparison, at $q=4$. The diagonal indicates unchanged phase.}
\label{fig:scc_visualization}
\end{figure}

The pointwise estimate stays within a factor $2.86$--$3.15$ of the nerve-map distance. At $q=4$, when landmarks merge, the maximum and RMS circular errors are $0.1073$ and $0.0654$ radians.

\section{Discussion and Future Work}
\label{Sec:Discussion}

Our main result, \Cref{Thm:Main-Theorem}, proves Lipschitz stability of the weighted nerve map under a uniform nondegeneracy assumption, with estimates in both metrics $\dTr$ and $\Wass$. Moreover, the experiments of \Cref{Sec:Experiments,Sec:Applications} suggest that the pointwise estimate $\Delta_L^r$ is the sharpest in practice, particularly for nonuniform weights, and each of the sharper two estimates in \Cref{Thm:Main-Theorem} can also be cheaper to compute exactly for large datasets.

We close with several potential directions for future work. First, although $\Delta_L^r$ need not satisfy the triangle inequality, one can take the infimum over finite chains to obtain a pseudometric $d_{\mathrm{ch}}^r$ satisfying $d_{\mathrm{ch}}^r\leq\Delta_L^r$ (\Cref{Appendix:Fail}). Whether this pseudometric additionally separates points under reasonable hypotheses remains open, and computing it may be difficult in practice. Another possibility is \emph{unbalanced optimal transport} \cite{LieroMielkeSavare2018, ChizatPeyreSchmitzerVialard2018, FigalliGigli2010, PiccoliRossi2014}, in which one is allowed to compare measures with different total mass by paying a cost for creating or destroying mass. Since our cost function $c_m$ admits an interpretation as a transport distance within $B(m,r)$ where an auxiliary point is allowed to absorb extra mass (\Cref{Lem:Cost-Pseudo}), unbalanced optimal transport may offer another route to defining a metric with comparable sharpness to $\Delta_L^r$.

Second, an important practical question is how to choose the scale $r$, landmarks $X$, and weights $\omega_X$ to obtain a well-conditioned partition of unity. For fixed $X$ and $r$, one can maximize the normalization factor $\inf \{ S_X^r(m) \mid m \in M\}$ over the weights, as in \Cref{Sec:Bio-Data}. However, increasing $r$ improves coverage and the normalization terms at the cost of losing local geometric information, while using more landmarks can improve coverage at greater computational cost. Optimizing the weights may also distort the distribution represented by the data. Finding selection procedures that balance these competing considerations remains an open problem.

Third, end-to-end stability for Laplacian eigenvector cascading \cite{MikePerea2019LaplacianEigenvectorCascading} and $K(G,n)$-valued coordinates remains open. Our estimates control the contribution of the nerve map, but in practice, the other data, e.g., eigenvectors, cocycle representatives, and smoothing also change with the landmarks. For cascading, spectral perturbation results such as Davis--Kahan \cite{davis1970rotation,YuWangSamworth2015} could help control changes in the relevant spectral subspaces. For $K(G,n)$ coordinates, induced matchings of barcodes \cite{BauerLesnick2015} and perturbation bounds for pseudoinverses \cite{Stewart1977} could be used to analyze persistent class selection and harmonic smoothing, respectively.

Finally, we would like probabilistic bounds for nerve maps built from random data, in a way that quantifies the dependence of stability on sample size, landmark selection, weights, and scale. Combining Wasserstein concentration estimates \cite{FournierGuillin2015} with bounds on expected cover radii \cite{ChazalGlisseLabruereMichel2015} could control both the transport error and the nondegeneracy constant in \Cref{Thm:Main-Theorem} with high probability.

\section*{Acknowledgements}
The authors thank Matan Harel for helpful discussions. We also thank Deane Yang for supplying a proof of \Cref{Prop:Ahlfors-Riemannian} on Mathematics Stack Exchange \cite{MSEPost}.

\bibliography{bibliography}
\bibliographystyle{plainnat}

\appendix
\crefalias{section}{appendix}

\section{Failure of the Triangle Inequality}
\label{Appendix:Fail}

The following example shows that the pointwise estimate $\Delta_L^r$ in \Cref{Thm:Main-Theorem} need not satisfy the triangle inequality, even when all measures have the same support.

\begin{example}\label{Ex:Delta-Triangle}
Let $M := \{0,1,3\} \subseteq \R$ with the Euclidean metric, and fix $r =2$, as in \Cref{Ex:Wass-Infty-Triangle}. Define the probability measures
\[
\begin{aligned}
\mu &= \tfrac18(\delta_0+6\delta_1+\delta_3),\\
\nu &= \tfrac18(3\delta_0+2\delta_1+3\delta_3),\\
\lambda &= \tfrac18(2\delta_0+5\delta_1+\delta_3).
\end{aligned}
\]
At $m=0,1$, the supports of the nerve maps are $\{0,1\}$, so $D_m=1$. The corresponding cost matrices, indexed by $0,1,3$, are
\[
\widetilde c_0=
\begin{pmatrix}
0&2&2\\
2&0&1\\
2&1&0
\end{pmatrix},
\qquad
\widetilde c_1=
\begin{pmatrix}
0&2&1\\
2&0&2\\
1&2&0
\end{pmatrix}.
\]
At $m=3$, we have $D_m=0$, so $\widetilde c_3=0$. Computing the transportation costs and normalization terms gives
\[
\Delta_L^2(\mu,\nu)=\frac34,\qquad
\Delta_L^2(\mu,\lambda)=\frac29,\qquad
\Delta_L^2(\lambda,\nu)=\frac12.
\]
Consequently,
\[
\Delta_L^2(\mu,\nu)=\frac34
>\frac29+\frac12
=\Delta_L^2(\mu,\lambda)+\Delta_L^2(\lambda,\nu).
\]
\end{example}

However, one can obtain a pseudometric from $\Delta_L^r$ using a standard construction from metric geometry that produces a pseudometric from any nonnegative, reflexive, and symmetric function. The construction works by taking the infimum of its sums along finite chains \cite[Chapter~I.1, 1.24]{bridsonHaefliger1999Metric}. In our setting, for a class $\mathcal A\subseteq\mathcal P(M)$, one may define
\[
d_{\mathrm{ch}}^r(\mu,\nu)
:=
\inf\left\{
\sum_{i=0}^{N-1}
\Delta_L^r(\mu_i,\mu_{i+1})
:
\begin{array}{l}
N\geq 1,\ \mu_i\in\mathcal A,\\
\mu_0=\mu,\ \mu_N=\nu
\end{array}
\right\}.
\]
By concatenating chains, $d_{\mathrm{ch}}^r$ is a pseudometric satisfying $d_{\mathrm{ch}}^r\leq\Delta_L^r$. However, its computation requires optimizing over arbitrarily many intermediate measures, making it difficult to use in practice. Moreover, on unrestricted classes of measures, separation can fail when the domains of two nerve maps are disjoint. Whether this construction separates points under suitable additional hypotheses is left for future work.

\section{Comparison to Wasserstein Transform}
\label{Appendix:WT}

The Wasserstein Transform of \cite{memoliWassersteinTransformArxiv2018} updates the metric on a dataset using \emph{localization operators}, which assign to each point a probability measure describing its neighborhood. Given $\alpha\in\mathcal P(M)$ and a continuous kernel $K:[0,\infty)\to[0,\infty)$, the \emph{kernel localization} of $\alpha$ at $m\in M$ is
\[
\alpha_m^K:=\frac{K(d_M(m,\cdot))}{S_\alpha^K(m)}\,\alpha,
\qquad
S_\alpha^K(m):=\int_M K(d_M(m,z))\,d\alpha(z),
\]
defined whenever $S_\alpha^K(m)>0$ \cite[Section~4.5.1]{memoliWassersteinTransformArxiv2018}. For $K_r(t):=|r-t|_+$ and $\alpha=\omega_X$, this is exactly our construction: $\alpha_m^{K_r}=\varphi_X^r(m)$ for $m\in X^r$.

For an $L$-Lipschitz kernel $K$, the proof of \cite[Lemma~48]{memoliWassersteinTransformArxiv2018} gives the following bound whenever both normalizations are positive; the same conclusion also follows by applying the lemma to $K+\varepsilon$ and letting $\varepsilon\downarrow0$:
\[
\Wass(\alpha_m^K,\beta_m^K)
\leq
\frac{2LD+\sup_{[0,D]}K}{\max\{S_\alpha^K(m),S_\beta^K(m)\}}\,
\Wass(\alpha,\beta)
\]
Moreover, \cite[Theorem~27]{memoliWassersteinTransformArxiv2018} obtains a uniform bound for strictly positive kernels using $\min_{[0,D]}K>0$. Our kernel $K_r$ is $1$-Lipschitz and bounded by $r$, but vanishes on $[r,\infty)$, so for $r\leq D$ the latter does not apply. The former gives, for $m\in X^r\cap Y^r$,
\[
\Wass\bigl(\varphi_X^r(m),\varphi_Y^r(m)\bigr)
\leq
\frac{2D+r}{\max\{S_X^r(m),S_Y^r(m)\}}\,
\Wass(\omega_X,\omega_Y).
\]

Our argument yields the same pointwise bound with $2D+r$ replaced by $3r$: by the proof of \Cref{Thm:Main-Theorem} and \Cref{Lem:Cost-Pseudo}, $\tilde c_m\leq c_m\leq 3r\,d_M$, hence
\[
\Wass\bigl(\varphi_X^r(m),\varphi_Y^r(m)\bigr)
\leq
\frac{T_{\tilde c_m}(\omega_X,\omega_Y)}{\max\{S_X^r(m),S_Y^r(m)\}}
\leq
\frac{3r}{\max\{S_X^r(m),S_Y^r(m)\}}\,
\Wass(\omega_X,\omega_Y).
\]
Thus our estimate removes the dependence on the diameter, improving the constant by the factor $(2D+r)/(3r)$. The reason is that our argument moves residual mass only within $\supp\varphi_X^r(m)\cup\supp\varphi_Y^r(m)\subseteq B(m,r)$, a set of diameter at most $2r$.

Refining the Wasserstein distance itself gives further improvements. The localized transport estimate replaces $3r\,\Wass(\omega_X,\omega_Y)$ with $\dTr(\omega_X,\omega_Y)$, a metric generating the same topology as $\Wass$ (\Cref{Prop:Localized-Metric-Properties}); its bound is attained in \Cref{Ex:Estimate-Sharp}. The pointwise estimate is sharper still, retaining the pointwise normalization and actual support diameter $D_m$, as illustrated in \Cref{Sec:Experiments,Sec:Applications}.

\end{document}